\documentclass[a4paper,11pt]{article}
\usepackage[T1]{fontenc}
\usepackage{babel}

\usepackage{authblk}
\usepackage[numbers,sort]{natbib} 

\usepackage{hyperref}
\usepackage{float}
\usepackage{amssymb,amscd,amsfonts,amsthm}
\usepackage[centertags]{amsmath}
\usepackage{enumerate}
\usepackage{latexsym,mathtools}
\usepackage{multicol}
\usepackage{tikz}

\newtheorem{theorem}{Theorem}

\newtheorem{corollary}[theorem]{Corollary}
\newtheorem{definition}[theorem]{Definition}
\newtheorem{example}[theorem]{Example}
\newtheorem{lemma}[theorem]{Lemma}
\newtheorem{proposition}[theorem]{Proposition}
\newtheorem{remark}[theorem]{Remark}

\newenvironment{claim}{\noindent\textbf{Claim.}}{\medskip}
\newenvironment{proofclaim}{\noindent{\emph{Proof.}}}{\hfill$(\square)$\medskip}

\newcommand{\Kneser}[2]{\ensuremath{K}({#1,#2})}

\newcommand{\neigho}[2]{\ensuremath{N_{#2}\left({#1}\right)}}
\newcommand{\neighc}[2]{\ensuremath{N_{#2}\left[{#1}\right]}}

\newcommand{\mindeg}[1]{\ensuremath{\delta\left({#1}\right)}}

\newcommand{\N}{\mathbb{N}}

\newcommand{\dom}[1]{\gamma\left(#1\right)}
\newcommand{\kdom}[2]{\gamma_{\times #2}\left(#1\right)}

\newcommand{\packing}[1]{\rho\left(#1\right)}

\usepackage{tabularx,environ}

\title{$k$-tuple domination on Kneser graphs
\footnote{Partially supported by PIP CONICET 1900, PICT-2020-03032 and PID 80020210300068UR.} 
}

\author[1,2]{Maria Gracia Cornet \footnote{e-mail: mgracia@fceia.unr.edu.ar}}
\author[1,2]{Pablo Torres \footnote{e-mail: ptorres@fceia.unr.edu.ar}}
\date{}

\affil[1]{\footnotesize Depto. de Matem\'atica, Facultad de Ciencias Exactas, Ing. y Agrimensura, Universidad Nacional de Rosario}
\affil[2]{\footnotesize Consejo Nacional de Investigaciones Científicas y Técnicas, Argentina}

\begin{document}
\maketitle
\sloppy

\begin{abstract}
This paper considers multiple domination on Kneser graphs.
We focus on $k$-tuple dominating sets, $2$-packings and the associated graph parameters $k$-tuple domination number and $2$-packing number.
In particular, we compute the $2$-packing number of Kneser graphs $\Kneser{3r-2}{r}$ and in odd graphs we obtain minimum $k$-tuple dominating sets of $K(7,3)$ and $K(11,5)$ for every $k$.
Besides, we determine the Kneser graphs $\Kneser{n}{r}$ with $k$-tuple domination number exactly $k+r$ and find all the minimum $k$-tuple dominating sets for these graphs, which generalize results for domination on Kneser graphs. Finally, we give a characterization of the $k$-tuple dominating sets of $\Kneser{n}{2}$ in terms of the occurrences of the elements in $[n]$, which allows us to obtain minimum sized $k$-tuple dominating sets of $\Kneser{n}{2}$ for $n\geq \Omega(\sqrt{k})$.
\end{abstract}

\section{Introduction}

Given a simple graph $G$,
let $\neigho{v}{G}$ denote the \emph{open neighbourhood} of a vertex $v$ in $G$ and $\neighc{v}{G}=\neigho{v}{G}\cup\{v\}$ the \emph{closed neighbourhood} of $v$ in $G$. 
When the graph $G$ is clear from the context, we may omit the subscripts and simply write $\neigho{v}{}$ and $\neighc{v}{}$. 
Furthermore, let $\mindeg{G}$ be the \emph{minimum degree} among all the vertices of graph $G$. Other definitions not given here can be
found in standard graph theory books such as \cite{west2001introduction}.
A \emph{dominating set} in $G$ is a subset $D\subseteq V(G)$ such that every vertex $v\in V(G)$ verifies that $|N[v]\cap D|\geq1$. The \emph{domination number} of $G$, denoted $\dom{G}$, is the minimum cardinality of a dominating set in $G$. Domination in graphs has been extensively studied in graph theory, and there is rich literature on this subject (see e.g. \cite{haynes2020topics,hhs-1998a,hhs-1998b}).

Some of the most studied variations of domination introduced a variable $k\in\N$, such as $k$-tuple domination.
Formally, given a graph $G$ and $k\in\N$, a set $D\subseteq V(G)$ is called a \emph{$k$-tuple dominating set} of $G$ if for every vertex $v\in V(G)$, we have $|\neighc{v}{G}\cap D|\geq k$. The \emph{$k$-tuple domination number} of $G$ is the minimum cardinality of a $k$-tuple dominating set of $G$, and is denoted by $\kdom{G}{k}$. 
A \emph{$\gamma_{\times k}$-set} is a $k$-tuple dominating set with cardinality $\kdom{G}{k}$. The $k$-tuple domination number is only defined for graphs with $\mindeg{G}\geq k-1$.
An excellent brief survey on $k$-tuple domination appears in the book Topics in Domination in Graphs \cite{haynes2020topics}, as part of a chapter devoted to the study of multiple domination.

Several bounds for the $k$-tuple domination number of graphs are given in \cite{martinez2022note,RAUTENBACH200798,zverovich2008k,chang2008upper,gagarin2008generalised,przybylo2013upper}. Liao and Chang \cite{liao2003k} studied the problem from an algorithmic point of view and proved that, fixed $k$, determining the $k$-tuple domination number is NP-complete for split graphs and for bipartite graphs.  Considering these unfavourable outcomes, exploring how this parameter behaves within classes of graphs exhibiting a nice combinatorial structure, like Kneser graphs, represents a natural direction for further research. In addition, another goal of this paper is to contribute to a deeper understanding on Kneser graphs. For instance, in Section \ref{packing} we show that $k$-tuple dominating sets in Kneser graphs are closely related to interesting combinatorial incidence structures like $2$-packings \cite{MeirMoon1975} in graphs, codes \cite{biggs1979} in graphs and intersecting set systems (e.g., Fano plane and Steiner systems).

For positive integers $r,n$ with $n\geq 2r$, the \emph{Kneser graph} $\Kneser{n}{r}$ has as vertices the $r$-subsets of an $n$-set and two vertices are adjacent in $\Kneser{n}{r}$ if and only if they are disjoint. This class of graphs gained prominence due to the Erdős-Ko-Rado theorem \cite{ekr-61}, which determined the independence number $\alpha(\Kneser{n}{r})$ of the Kneser graph $\Kneser{n}{r}$ to be $\binom{n-1}{r-1}$, as a result on extremal combinatorics. Lovász's proof of Kneser's conjecture \cite{lov-78,kneser1955aufgabe}, later complemented by Matoušek's combinatorial proof \cite{mat-04}, provided a determination of the chromatic number of Kneser graphs. Numerous other graph invariants have been investigated on Kneser graphs by various authors \cite{TreewidthKneser,qAnalog,GenPosKneser,bedo2023geodetic}.
One of them is the domination number \cite{hw-2003,ivanvco1993domination,ostergaard2014bounds,gorodezky2007dominating}. In particular, the value of the domination number of $\Kneser{n}{r}$ is determined for $n$ large enough in {Theorem \ref{TheoDomBigN}}.

\begin{figure}[ht]
    \centering
    \begin{tikzpicture}
            \foreach \x in {0,1,...,4}{
            \draw ({90+\x*72}:1) -- ({90+\x*72}:2);
            \draw ({90+\x*72}:2) -- ({90+(\x+1)*72}:2);
            \draw ({90+\x*72}:1) -- ({90+(\x+2)*72}:1);
            \filldraw ({90+\x*72}:1) circle (2pt);
            \filldraw ({90+\x*72}:2) circle (2pt);
            }

            \foreach \x in {0,1,...,4}{
            \begin{scope}[shift={({90+\x*72}:1)}]
                \filldraw[fill=white, thick] (0,0) circle (0.2);
                \foreach \y in {0,1,...,4}{
                \draw (0,0) -- ({90-\y*72}:0.2);
                }
            \end{scope}
            \begin{scope}[shift={({90+\x*72}:2)}]
                \filldraw[fill=white, thick] (0,0) circle (0.2);
                \foreach \y in {0,1,...,4}{
                \draw (0,0) -- ({90-\y*72}:0.2);
                }
            \end{scope}
            }

            \begin{scope}[shift={(90:2)}]
                \filldraw[fill=gray!70] (0,0) -- (90:0.2) arc (90:18:0.2) -- cycle;
                \filldraw[fill=gray!70] (0,0) -- (18:0.2) arc (18:-54:0.2) -- cycle;
            \end{scope}
            
            \begin{scope}[shift={(18:2)}]
                \filldraw[fill=gray!70] (0,0) -- (-126:0.2) arc (-126:-198:0.2) -- cycle;
                \filldraw[fill=gray!70] (0,0) -- (-198:0.2) arc (-198:-270:0.2) -- cycle;
            \end{scope}
            
            \begin{scope}[shift={(-54:2)}]
                \filldraw[fill=gray!70] (0,0) -- (18:0.2) arc (18:-54:0.2) -- cycle;
                \filldraw[fill=gray!70] (0,0) -- (-54:0.2) arc (-54:-126:0.2) -- cycle;
            \end{scope}

            \begin{scope}[shift={(-126:2)}]
                \filldraw[fill=gray!70] (0,0) -- (90:0.2) arc (90:18:0.2) -- cycle;
                \filldraw[fill=gray!70] (0,0) -- (-198:0.2) arc (-198:-270:0.2) -- cycle;
            \end{scope}

            \begin{scope}[shift={(-198:2)}]
                \filldraw[fill=gray!70] (0,0) -- (-54:0.2) arc (-54:-126:0.2) -- cycle;
                \filldraw[fill=gray!70] (0,0) -- (-126:0.2) arc (-126:-198:0.2) -- cycle;
            \end{scope}
            
            \begin{scope}[shift={(90:1)}]
                \filldraw[fill=gray!70] (0,0) -- (-54:0.2) arc (-54:-126:0.2) -- cycle;
                \filldraw[fill=gray!70] (0,0) -- (-198:0.2) arc (-198:-270:0.2) -- cycle;
            \end{scope}
            
            \begin{scope}[shift={(18:1)}]
                \filldraw[fill=gray!70] (0,0) -- (90:0.2) arc (90:18:0.2) -- cycle;
                \filldraw[fill=gray!70] (0,0) -- (-54:0.2) arc (-54:-126:0.2) -- cycle;
            \end{scope}
            
            \begin{scope}[shift={(-54:1)}]
                \filldraw[fill=gray!70] (0,0) -- (90:0.2) arc (90:18:0.2) -- cycle;
                \filldraw[fill=gray!70] (0,0) -- (-126:0.2) arc (-126:-198:0.2) -- cycle;
            \end{scope}
            
            \begin{scope}[shift={(-126:1)}]
                \filldraw[fill=gray!70] (0,0) -- (18:0.2) arc (18:-54:0.2) -- cycle;
                \filldraw[fill=gray!70] (0,0) -- (-126:0.2) arc (-126:-198:0.2) -- cycle;
            \end{scope}
            
            \begin{scope}[shift={(-198:1)}]
                \filldraw[fill=gray!70] (0,0) -- (18:0.2) arc (18:-54:0.2) -- cycle;
                \filldraw[fill=gray!70] (0,0) -- (-198:0.2) arc (-198:-270:0.2) -- cycle;
            \end{scope}

        \end{tikzpicture}
    
    \caption{Kneser graph $\Kneser{5}{2}$}
    \label{fig: Petersen}
\end{figure}
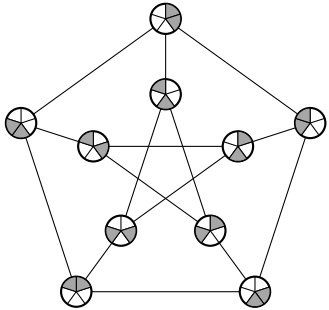

\begin{theorem}\cite[Theorem 2.2]{ostergaard2014bounds}\label{TheoDomBigN}
    If $n\geq r^2+r$, then $\dom{\Kneser{n}{r}}{}=r+1$.
\end{theorem}

It is also known that $\dom{\Kneser{n}{2}}{}=3$ for $n\geq4$ \cite{ivanvco1993domination}. Furthermore, for the remaining cases of Kneser graphs partial results are provided, but a complete solution for the domination number, similar to the case of chromatic and independence numbers, has not yet been achieved.

After these results, different variations of domination have been explored on Kneser graphs. For instance, in 2019 Bre\v{s}ar et al. \cite{bresar2019grundy} studied four different types of Grundy domination numbers and the related zero forcing numbers on Kneser graphs. In 2023, P. Jalilolghadr and A. Behtoei computed the total dominator chromatic number of the Kneser graph $\Kneser{n}{2}$ for each $n\geq5$. Also in 2023, three Roman domination graph invariants on Kneser graphs are studied by T. Zec and M. Grbi\'c \cite{zec2023several}. They obtained the exact values of the mentioned invariants for $\Kneser{n}{2}$ and also for $\Kneser{n}{r}$ if $n$ is large enough with respect to $r$. 
A closely related concept to $k$-tuple domination is \emph{$k$-domination}.
Given a graph $G$ and $k\in\mathbb{N}$, a set $D\subseteq V(G)$ is a  \emph{$k$-dominating set} if for every vertex $u\notin D$ it holds $|N[u]\cap D|\geq k$. The \emph{$k$-domination number} of a graph $G$ is the minimum cardinality of a $k$-dominating set of $G$. Results on $k$-domination for Kneser graphs appear in \cite{cornet2023kdom} after a first version of the present paper \cite{ktuple2023v1} was released.

As we have mentioned, in this work we focus on the study of $k$-tuple domination on Kneser graphs. In Section \ref{packing} we analyze $k$-tuple dominating sets of Kneser graphs, from its relationship with $2$-packings in graphs. 
We show that $2$-packings on Kneser graphs are closely connected to intersecting set families. Using this link, we compute the $2$-packing number of $\Kneser{n}{r}$ for $n=3r-2$ and $r\geq3$. 
Then we obtain $\gamma_{\times k}$-sets for the odd graphs $\Kneser{7}{3}$ and $\Kneser{11}{5}$ for every $k$, from the Steiner systems $S(2,3,7)$ (Fano plane) and $S(4,5,11)$ respectively. Furthermore, in the case of $\Kneser{7}{3}$, we give all the $\gamma_{\times k}$-sets.

In Section \ref{sec:general} we generalize the result in Theorem \ref{TheoDomBigN} (notice that $\kdom{G}{1}=\dom{G}$) by proving that $\kdom{\Kneser{n}{r}}{k}=k+r$ if and only if $n\geq r(k+r)$ and we characterize the $\gamma_{\times k}$-sets for these cases. Besides, we prove that $\kdom{\Kneser{n}{r}}{k}$ is not decreasing with respect to $n$, so we conclude that $\kdom{\Kneser{n}{r}}{k}$ is at least $k+r+1$ if $n< r(k+r)$. In addition, we calculate $\kdom{\Kneser{n}{r}}{k}$ for $k$ large enough.

Finally, in Section \ref{secKn2} we give a characterization of the $k$-tuple dominating sets for $\Kneser{n}{2}$ that lead us to provide $\kdom{\Kneser{n}{2}}{k}$ for $k\leq a\,\frac{n-3}{4}$, with $a=n-4$ if $n$ is even and $a=n-6+(n\ mod\ 4)$ if $n$ is odd.
We also design an ILP formulation to obtain some $k$-tuple dominating sets whose cardinality meet the lower bounds given previously. Table \ref{tab: tabla K(n,2)} shows the results for $n\leq26$ and $k\leq 60$.

\section{\texorpdfstring{$2$-packings and $k$-tuple dominating sets on odd graphs}{2-packings and k-tuple dominating sets of odd graphs}}\label{packing}

Packings and dominating sets have been extensively studied in graphs. In this regard, $2$-packings have been received particular attention and determining the $2$-packing number of classes of graphs is a challenging and interesting topic, see e.g. \cite{2packing-JSR-2012,2packing-BKPG-2022,2packingCG-2017,2packingCKMR-2011}. Given a graph $G$, a subset $S\subseteq V(G)$ is 
called a \emph{$2$-packing} of $G$ if for every pair of vertices $u,v\in S$, their closed neighborhoods satisfy $\neighc{u}{}\cap \neighc{v}{}=\emptyset$. In other words, no two vertices in the $2$-packing have any common neighbor. The \emph{$2$-packing number} of a graph $G$, denoted by $\packing{G}$, is the maximum cardinality of a $2$-packing in $G$.  In \cite{martinez2022note}, the author establishes a relationship between the $2$-packing number and the $k$-tuple domination number.

\begin{theorem}\cite[Theorem 2.3]{martinez2022note}
\label{thm: cota packing}
    Let $k\geq 2$ be an integer. For any graph $G$ of order $n$ and $\mindeg{G}\geq k$,
    $$k\packing{G}\leq \kdom{G}{k}\leq n-\packing{G}.$$
\end{theorem}

It is known that when $n\geq 3r-1$ with $r\geq 2$, the diameter of the Kneser graph $\Kneser{n}{r}$ is equal to $2$ \cite{valencia2005diameter}.
Therefore, in these cases the $2$-packing number of the Kneser graph $\Kneser{n}{r}$ is equal to $1$.
Thus, for the remaining of this section, we consider $2r+1\leq n\leq 3r-2$.
Similarly to other works on domination on Kneser graphs, we use $\packing{n,r}$ to denote the $2$-packing number of the Kneser graph $\Kneser{n}{r}$, and $\kdom{n,r}{k}$ to denote the $k$-tuple domination number of the Kneser graph $\Kneser{n}{r}$. For positive integers $r\leq n$, we denote by $[r..n]$ and $[n]$ the sets $\{r,\ldots,n\}$ and $\{1,\ldots,n\}$, respectively.

Let $S$ be a $2$-packing in $\Kneser{n}{r}$. If $u$ and $v$ are two vertices in $S$ then $u\cap v\neq \emptyset$. Besides, notice that if $|u\cap v|>(3r-1)-n$, then
$$|[n]\setminus(u\cup v)|=|[n]|-|u|-|v|+|u\cap v|\geq r,$$
This implies that there exists a vertex $w\in \neigho{u}{}\cap\neigho{v}{}$, contradicting the fact that $S$ is a $2$-packing.

\begin{remark}\label{rem: intersecciones en un packing}
    Let $2r+1\leq n\leq 3r-2$. A set $S$ of $r$-subsets of $[n]$ is a $2$-packing of $\Kneser{n}{r}$ if and only if for every pair $u,v\in S$, it holds that
    $$1\leq|u\cap v|\leq (3r-1)-n.$$
\end{remark}

Some results on extremal combinatorics and intersecting families can be used to give upper bounds for the $2$-packing number of Kneser graphs. In fact, from \cite{ramanan1997proof}, we have:
$$\packing{n,r}\leq \binom{n}{(3r-1)-n}$$
and
$$\packing{n,r}\leq \sum_{i=0}^{(3r-1)-n} \binom{n-1}{i}.$$
In the case $n=3r-2$ both bounds remain
$$\packing{n,r}\leq n,$$
which is tight only for $r=3$, as we will show in Theorem \ref{teopacking}.

We introduce the following notation, which will be used throughout the paper. Given a set of vertices $D$ in $\Kneser{n}{r}$ and $x\in [n]$, the occurrences of the element $x$ in $D$, denoted by $i_x(D)$, represent the number of vertices in $D$ that contain the element $x$. In other words, $i_x(D)$ is the cardinality of the set $\{u\in D : x\in u\}$. For a positive integer $a$, we define $X_a(D)$ as the set of elements in $[n]$ such that their occurrences in $D$ are equal to $a$, i.e., $X_a(D)=\{x\in [n] : i_x(D)=a\}$. Similarly, we define $X_a^\geq(D)=\{x\in [n] : i_x(D)\geq a\}$, and $X_a^\leq(D)=\{x\in [n] : i_x(D)\leq a\}$. When the set $D$ is clear from the context, we shall omit it in the notation. It is important to note that the sum of the occurrences of all elements in $D$ is equal to $r$ times the cardinality of $D$, i.e., $\sum_{x\in [n]} i_x=r|D|$.

\begin{theorem}\label{teopacking}
    If $n=3r-2$ and $r\geq 3$, then
    \begin{equation}
        \packing{n,r}=\begin{cases}
        7, & \text{if }r=3,\\
        5, & \text{if }r=4,\\
        3, & \text{if }r\geq 5.
        \end{cases}
    \end{equation}
    
    \begin{proof}
    
    Let $r\geq 3$, $n=3r-2$ and let $S$ be a $2$-packing of $\Kneser{n}{r}$. From Remark \ref{rem: intersecciones en un packing}, we have that for every pair of vertices $u,v\in S$, it holds $|u\cap v|=1$.
    
    Let us see that $i_x\leq 3$ for every $x\in[n]$. In fact, suppose $i_a\geq 4$ for some $a\in[n]$, and let $u_1,u_2,u_3,u_4$ be four vertices in $S$ such that $a\in u_j$ for $j\in [4]$. We have $|u_j\setminus \{a\}|=r-1$ and $\left(u_j\setminus \{a\}\right)\cap \left(u_\ell\setminus \{a\}\right)=\emptyset$ for $j,\ell\in [4]$ with $j\neq \ell$. Then
    $$n\geq \left| \bigcup_{j=1}^4 u_j \right|=1+\sum_{j=1}^4\left|u_j\setminus \{a\}\right|=4r-3=n+\underbrace{r-1}_{>0},$$
    which cannot be true. Thus, $i_x\leq 3$ for every $x\in[n]$.
    Then we have $|S|\leq \dfrac{3n}{r}$.

    If $r=3$, then $n=7$ and the cardinality of a $2$-packing $S$ is at most $7$. In fact, a set of $7$ $3$-subsets of $[7]$ which mutually intersect in exactly one element, is a Fano plane (see Figure \ref{fig: K(7,3)}). Thus, we have $\packing{7,3}=7$.

    \begin{figure}[ht]
        \centering
        \begin{tikzpicture}
            \node (145) at ({-90+0*(360/7)}:1.3) {};
            \node (147) at ({-90+1*(360/7)}:1.3) {};
            \node (347) at ({-90+2*(360/7)}:1.3) {};
            \node (367) at ({-90+3*(360/7)}:1.3) {};
            \node (236) at ({-90+4*(360/7)}:1.3) {};
            \node (256) at ({-90+5*(360/7)}:1.3) {};
            \node (125) at ({-90+6*(360/7)}:1.3) {};
            
            
            \node (136) at ({90+0*(360/7)}:2.4) {};
            \node (246) at ({90+1*(360/7)}:2.4) {};
            \node (257) at ({90+2*(360/7)}:2.4) {};
            \node (135) at ({90+3*(360/7)}:2.4) {};
            \node (146) at ({90+4*(360/7)}:2.4) {};
            \node (247) at ({90+5*(360/7)}:2.4) {};
            \node (357) at ({90+6*(360/7)}:2.4) {};
            
            \node (127) at ({90+0*(360/7)}:3.5) {};
            \node (345) at ({90+1*(360/7)}:3.5) {};
            \node (167) at ({90+2*(360/7)}:3.5) {};
            \node (234) at ({90+3*(360/7)}:3.5) {};
            \node (567) at ({90+4*(360/7)}:3.5) {};
            \node (123) at ({90+5*(360/7)}:3.5) {};
            \node (456) at ({90+6*(360/7)}:3.5) {};
            
            \node (245) at ({90+360/21+0*(360/7)}:4) {};
            \node (157) at ({90+360/21+1*(360/7)}:4) {};
            \node (134) at ({90+360/21+2*(360/7)}:4) {};
            \node (467) at ({90+360/21+3*(360/7)}:4) {};
            \node (237) at ({90+360/21+4*(360/7)}:4) {};
            \node (356) at ({90+360/21+5*(360/7)}:4) {};
            \node (126) at ({90+360/21+6*(360/7)}:4) {};
            
            \node (137) at ({90+2*360/21+0*(360/7)}:4) {};
            \node (346) at ({90+2*360/21+1*(360/7)}:4) {};
            \node (267) at ({90+2*360/21+2*(360/7)}:4) {};
            \node (235) at ({90+2*360/21+3*(360/7)}:4) {};
            \node (156) at ({90+2*360/21+4*(360/7)}:4) {};
            \node (124) at ({90+2*360/21+5*(360/7)}:4) {};
            \node (457) at ({90+2*360/21+6*(360/7)}:4) {};

            \draw (123) -- (456);
            \draw (123) -- (457);
            \draw (123) -- (467);
            \draw (123) -- (567);
            
            \draw (124) -- (356);
            \draw (124) -- (357);
            \draw (124) -- (567);
            
            \draw (125) -- (347);
            \draw (125) -- (367);
            
            \draw (126) -- (345);
            \draw (126) -- (457);
            \draw (126) -- (357);
            
            \draw (127) -- (456);
            \draw (127) -- (345);
            \draw (127) -- (346);
            \draw (127) -- (356);
            
            \draw (134) -- (257);
            \draw (134) -- (267);
            \draw (134) -- (567);
            
            \draw (135) -- (267);
            \draw (135) -- (467);
            
            \draw (136) -- (245);
            \draw (136) -- (457);
            
            \draw (137) -- (245);
            \draw (137) -- (246);
            \draw (137) -- (456);
            
            \draw (145) -- (236);
            \draw (145) -- (367);
            
            \draw (146) -- (235);
            \draw (146) -- (237);
            
            \draw (147) -- (236);
            \draw (147) -- (256);
            
            \draw (156) -- (234);
            \draw (156) -- (237);
            \draw (156) -- (247);
            
            \draw (157) -- (234);
            \draw (157) -- (246);
            \draw (157) -- (346);
            
            \draw (167) -- (234);
            \draw (167) -- (235);
            \draw (167) -- (245);
            \draw (167) -- (345);

            \draw (234) -- (567);
            \draw (235) -- (467);
            \draw (237) -- (456);
            \draw (247) -- (356);
            \draw (256) -- (347);
            \draw (257) -- (346);
            \draw (267) -- (345);

            \draw (136) to[bend left=25] (247);
            \draw (247) to[bend left=25] (135);
            \draw (135) to[bend left=25] (246);
            \draw (246) to[bend left=25] (357);
            \draw (357) to[bend left=25] (146);
            \draw (146) to[bend left=25] (257);
            \draw (257) to[bend left=25] (136);

            \draw (236) to[bend right =30] (457);
            \draw (256) to[bend right =30] (137);
            \draw (125) to[bend right =30] (346);
            \draw (145) to[bend right =30] (267);
            \draw (147) to[bend right =30] (235);
            \draw (347) to[bend right =30] (156);
            \draw (367) to[bend right =30] (124);
            
            \draw (236) to[bend left =30] (157);
            \draw (256) to[bend left =30] (134);
            \draw (125) to[bend left =30] (467);
            \draw (145) to[bend left =30] (237);
            \draw (147) to[bend left =30] (356);
            \draw (347) to[bend left =30] (126);
            \draw (367) to[bend left =30] (245);
            
            

            \foreach \x in {123,124,125,126,127,134,135,136,137,145,146,147,156,157,167,234,235,236,237,245,246,247,256,257,267,345,346,347,356,357,367,456,457,467,567}{

            \begin{scope}[shift={(\x)}]
                \filldraw[fill = white, thick] (0,0) circle (0.2);
                \foreach \y in {0,1,...,6}{
                \draw (0,0) -- ({90-360/7*\y}:0.2);
                }
                \filldraw[fill = gray!70] (0,0) -- ({90-360/7*(floor(\x/100)-1)}:0.2) arc({90-360/7*(floor(\x/100)-1)}:{90-360/7*floor(\x/100)}:0.2) -- cycle;
                \filldraw[fill = gray!70] (0,0) -- ({90-360/7*(floor(\x/10)-floor(\x/100)*10-1)}:0.2) arc({90-360/7*(floor(\x/10)-floor(\x/100)*10-1)}:{90-360/7*(floor(\x/10)-floor(\x/100)*10)}:0.2) -- cycle;
                \filldraw[fill = gray!70] (0,0) -- ({90-360/7*(\x-floor(\x/10)*10-1)}:0.2) arc({90-360/7*(\x-floor(\x/10)*10-1)}:{90-360/7*(\x-floor(\x/10)*10)}:0.2) -- cycle;
            \end{scope}}

            \foreach \x in {124,156,235,267,346,137,457}{

            \begin{scope}[shift={(\x)}]
                \filldraw[fill = white, thick] (0,0) circle (0.2);
                \foreach \y in {0,1,...,6}{
                \draw (0,0) -- ({90-360/7*\y}:0.2);
                }
                \filldraw[fill = red] (0,0) -- ({90-360/7*(floor(\x/100)-1)}:0.2) arc({90-360/7*(floor(\x/100)-1)}:{90-360/7*floor(\x/100)}:0.2) -- cycle;
                \filldraw[fill = red] (0,0) -- ({90-360/7*(floor(\x/10)-floor(\x/100)*10-1)}:0.2) arc({90-360/7*(floor(\x/10)-floor(\x/100)*10-1)}:{90-360/7*(floor(\x/10)-floor(\x/100)*10)}:0.2) -- cycle;
                \filldraw[fill = red] (0,0) -- ({90-360/7*(\x-floor(\x/10)*10-1)}:0.2) arc({90-360/7*(\x-floor(\x/10)*10-1)}:{90-360/7*(\x-floor(\x/10)*10)}:0.2) -- cycle;
            \end{scope}}

            \begin{scope}[shift={(7.5,-0.5)}]
                \draw (90:2.5) -- (210:2.5) -- (-30:2.5) -- cycle;
                \draw (90:2.5) -- (-90:1.25);
                \draw (-30:2.5) -- (150:1.25);
                \draw (210:2.5) -- (30:1.25);
                \draw (0,0) circle (1.25);
            
                \filldraw (0,0) circle (2.5pt) node[xshift=6,yshift=10] {$6$};
                \filldraw (90:2.5) circle (2.5pt) node[yshift=10] {$4$};
                \filldraw (210:2.5) circle (2.5pt) node[xshift=-6,yshift=-10] {$7$};
                \filldraw (-30:2.5) circle (2.5pt) node[xshift=6,yshift=-10] {$1$};
                \filldraw (-90:1.25) circle (2.5pt) node[yshift=-10] {$3$};
                \filldraw (30:1.25) circle (2.5pt) node[xshift=6,yshift=8] {$2$};
                \filldraw (150:1.25) circle (2.5pt) node[xshift=-6,yshift=8] {$5$};
            \end{scope}

        \end{tikzpicture}
        \caption{An maximum $2$-packing in $\Kneser{7}{3}$ and the corresponding Fano plane}
        \label{fig: K(7,3)}
    \end{figure}
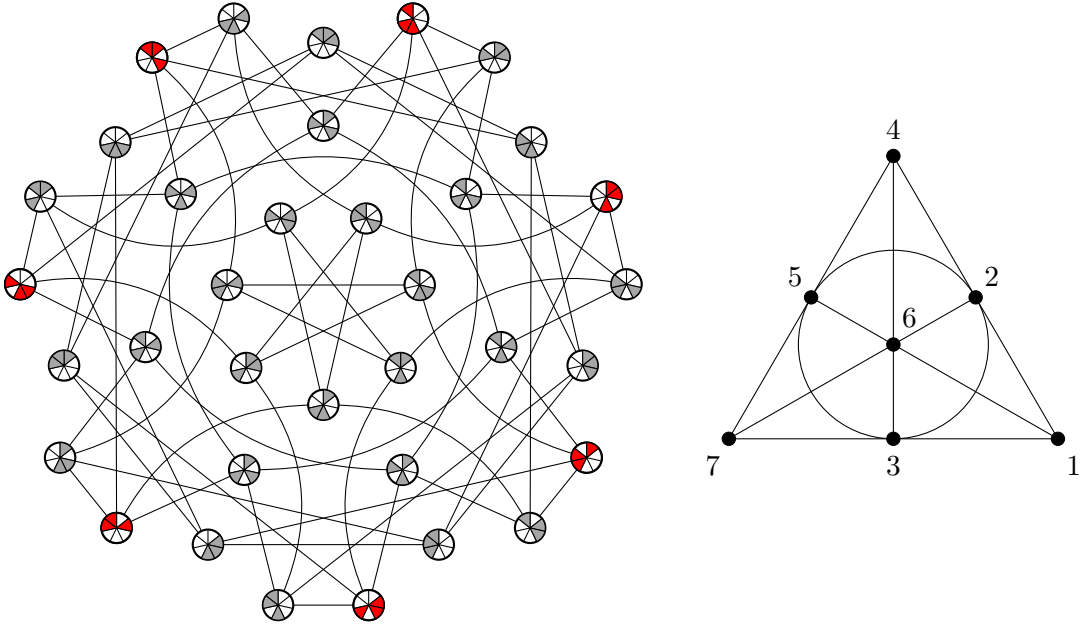
    
    If $r\geq 4$, let us suppose that $i_a=3$ for some $a\in[n]$ and let the vertices $u_1,u_2,u_3\in S$ be such that $a\in u_j$ for each $j$. Without loss of generality, let us consider $a=1$ and
    $$u_1=\{1\}\cup [2..r],\qquad u_2=\{1\}\cup [r+1..2r-1],\qquad u_3=\{1\}\cup [2r..3r-2].$$
    Let $w$ be another vertex in $S$. $1\notin w$ since $i_1$ is exactly $3$. Then $w=\{a_1,\ldots,a_{r}\}$ with $a_j\in[n]\setminus\{1\}$ for $j\in [r]$. Since $\left\{u_i\setminus{\{1\}}\right\}_{i=1}^3$ is a partition of $[n]\setminus{\{1\}}$, and $r\geq 4$, by the pigeonhole principle there exist at least two elements $a_j$ and $a_{\ell}$ that belongs to $u_k$ for some $k\in [3]$. This is, $|u_k\cap w|\geq 2$ but this leads to a contradiction with the fact that $S$ is a $2$-packing.
    
    Therefore, if there exists $a\in [n]$ such that $i_a=3$, then $|S|\leq 3$. On the contrary, if $i_a\leq 2$ for every $a\in [n]$, we have $|S|\leq \dfrac{2n}{r}$.
    
    If $r=4$, then $n=10$ and the cardinality of a $2$-packing $S$ is at most $5$. In fact, we have $\packing{10,4}=5$ and if $S$ is an maximum $2$-packing of $\Kneser{10}{4}$, then, up to automorphism, we have $S=\big\{ \{1,2,3,4\},\{1,5,6,7\},\{2,5,8,9\},\{3,6,8,10\},\{4,7,9,10\} \big\}$.
    
    If $r\geq 5$, suppose that $i_a\leq 2$ for every $a\in[n]$. Let $u_1,u_2,u_3\in S$. Without loss of generality, we may assume that
    $$u_1=\{1,2\}\cup [4..r+1],\qquad u_2=\{1,3\}\cup [r+2..2r-1],\qquad u_3=\{2,3\}\cup [2r..3r-3].$$
    A fourth vertex in $S$ must contain exactly one element from each set $u_1,u_2,u_3$ and two elements from $[n]\setminus\left(u_1\cup u_2\cup u_3\right)$. However,
    $$\left| [n]\setminus\left(u_1\cup u_2\cup u_3\right) \right|=(3r-2)-(3r-3)=1.$$ 
    We can conclude that $|S|\leq 3$, and $S=\{u_1,u_2,u_3\}$ is a $2$-packing of maximum cardinality. 
    
    Therefore, if $r\geq 5$, $\packing{3r-2,r}=3$.
    \end{proof}
\end{theorem}

As a by-product of the previous proof, we have the following result.

\begin{corollary}
    If $n=3r-2$ with $r\geq 5$, and $S$ is an maximum $2$-packing of $\Kneser{n}{r}$, then up to automorphism we have
    $$S=\big\{ \{1,2\}\cup [4..r+1],\{1,3\}\cup [r+2..2r-1],\{2,3\}\cup [2r..3r-3] \big\},$$
    or 
    $$S=\big\{ \{1\}\cup [2..r],\{1\}\cup [r+1..2r-1],\{1\}\cup [2r..3r-2] \big\}.$$
\end{corollary}

Let us notice that $2$-packings on Kneser graphs are related to Steiner systems. A
Steiner system $S(t,r,n)$ is a collection of $r$-subsets of $[n]$, called blocks, with the property that each $t$-subset of $[n]$ is contained in exactly one block. It is not hard to see that for $2r+1\leq n\leq3r-2$, if a Steiner system $S(3r-n,r,n)$ exists, then 
$$\packing{n,r}\leq|S(3r-n,r,n)|.$$

This relationship is useful to study $2$-packings in odd graphs. To this end, we consider perfect $1$-codes in graphs \cite{hammond1975perfect}. A subset of vertices $C$ of a graph $G$ is a perfect $1$-code of $G$ if the family of closed neighbourhood $\left\{N[v]\right\}_{v\in C}$ is a partition of $V(G)$. Notice that if a perfect $1$-code of a graph $G$ exists, then it is also an maximum $2$-packing of $G$.
In particular, it is known \cite{hammond1975perfect} that a set $C$ is a Steiner system $S(r-1,r,2r+1)$ if and only if $C$ is a perfect $1$-code in the odd graph $\Kneser{2r+1}{r}$. 
Thus, the case $r=3$ in Theorem $\ref{teopacking}$ also follows from the fact that the Fano plane is a Steiner system $S(2,3,7)$. Similarly, since a Steiner system $S(4,5,11)$ exists \cite{barrau1908combinatory}, then the odd graph $\Kneser{11}{5}$ has a perfect $1$-code and $$\packing{11,5}=66.$$

The well known conjecture due to Biggs \cite{biggs1979} asserts that there is no perfect $1$-code if $r\neq3,5$. 
Although it has been verified for some values of $r$, it has not yet been settled in general.
For an exhaustive study of Steiner systems we refer to \cite{colbourn2006steiner} and references therein.

Notice that the existence of a perfect $1$-code of the odd graph $\Kneser{2r+1}{r}$ gives the following lower bound for the $k$-tuple domination number
$$\kdom{2r+1,r}{k}\geq k\packing{2r+1,r}=\frac{\binom{2r+1}{r}\, k}{r+2}$$
which we will show is tight for $r=3,5$ and each $k$.

\begin{lemma} \label{lem: cota packing para vertice-transitivos}
    Let $G$ be a vertex-transitive graph. If $\kdom{G}{k}=k\packing{G}$ and $D$ is a $\gamma_{\times k}$-set, then it holds $|N[v]\cap D|=k$ for every $v\in V(G)$.

    \begin{proof}
        Let $S$ be an maximum $2$-packing and $D$ a $\gamma_{\times k}$-set of $G$. Notice that $|N[u]\cap D|\geq k$ for every vertex $u$ as $D$ is a $k$-tuple dominating set. On the other hand, the sets $\{N[u]\}_{u\in S}$ are pairwise disjoint. So, we have 
        \begin{equation}
            \kdom{G}{k}=|D|\geq \sum_{u\in S} \underbrace{|N[u]\cap D|}_{\geq k}\geq k|S|=k\packing{G}=\kdom{G}{k}.
        \end{equation}
        In consequence, $|N[u]\cap D|=k$ for every $u\in S$.

        Besides, let $v\in V(G)$. Since $G$ is a vertex-transitive graph, it is possible to find an maximum $2$-packing $S'$ such that $v\in S'$. Thus, $|N[v]\cap D|=k$.
    \end{proof}
\end{lemma}

As it is well-known that Kneser graphs are vertex-transitive graphs, we have that if $\kdom{n,r}{k}=k\packing{n,r}$ and $D$ is a $\gamma_{\times k}$-set of $\Kneser{n}{r}$, then $|N[u]\cap D|=k$ for every $u\in V(\Kneser{n}{r})$.

\begin{remark}
\label{rem: complemento de kupla dom cota packing}
    Let us notice that if $D_k$ is a $k$-tuple dominating set of $\Kneser{n}{r}$ with cardinality $|D_k|=k\packing{n,r}$, then the set $V(\Kneser{n}{r})\setminus D_k$ is a $\tilde{k}$-tuple dominating set of $\Kneser{n}{r}$, with $\tilde{k}=\binom{n-r}{r}+1-k$. In fact, for every vertex $u\in V(\Kneser{n}{r})$, by Lemma \ref{lem: cota packing para vertice-transitivos} we have
    $$|N[u]\cap (V(\Kneser{n}{r})\setminus D_k)|=\underbrace{|N[u]|}_{\binom{n-r}{r}+1}-\underbrace{|N[u]\cap D_k|}_{k}=\binom{n-r}{r}+1-k.$$
\end{remark}

It is known that there exist exactly two disjoint Fano planes. Thus the union of these two Fano plane is a $\gamma_{\times 2}$-set of $\Kneser{7}{3}$. Moreover, we provide $\gamma_{\times k}$-sets of $\Kneser{7}{3}$ for each $k\in[5]$.

\begin{remark}\label{remarkkdom73}
    For $r=3$ and $n=7$, from Theorem \ref{teopacking} it turns out that $\packing{7,3}=7$. A bound for $\kdom{7,3}{k}$ is given by Theorem \ref{thm: cota packing} and we have $\kdom{7,3}{k}\geq 7k$. In fact, for these values of $n$ and $r$ we shall see that this bound is tight and we have $\kdom{7,3}{k}=7k$ for $1\leq k\leq 5$. 

    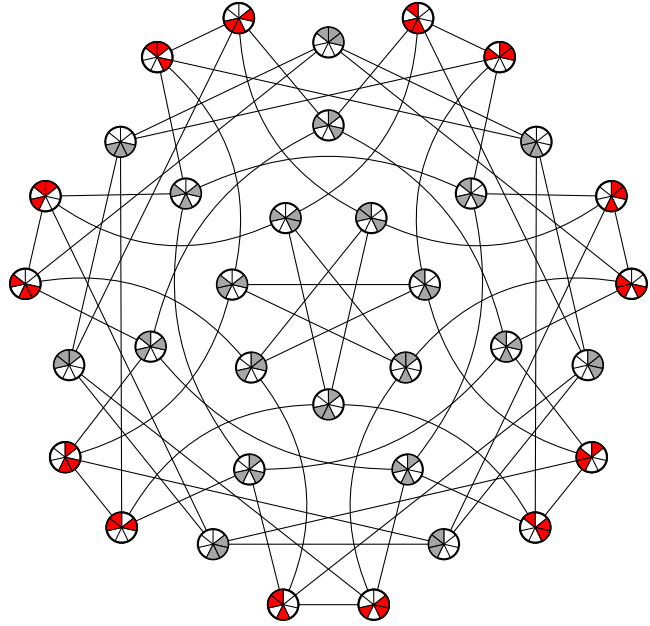
\begin{figure}[H]
        \centering
        \begin{tikzpicture}
            \node (145) at ({-90+0*(360/7)}:1.3) {};
            \node (147) at ({-90+1*(360/7)}:1.3) {};
            \node (347) at ({-90+2*(360/7)}:1.3) {};
            \node (367) at ({-90+3*(360/7)}:1.3) {};
            \node (236) at ({-90+4*(360/7)}:1.3) {};
            \node (256) at ({-90+5*(360/7)}:1.3) {};
            \node (125) at ({-90+6*(360/7)}:1.3) {};
            
            
            \node (136) at ({90+0*(360/7)}:2.4) {};
            \node (246) at ({90+1*(360/7)}:2.4) {};
            \node (257) at ({90+2*(360/7)}:2.4) {};
            \node (135) at ({90+3*(360/7)}:2.4) {};
            \node (146) at ({90+4*(360/7)}:2.4) {};
            \node (247) at ({90+5*(360/7)}:2.4) {};
            \node (357) at ({90+6*(360/7)}:2.4) {};
            
            \node (127) at ({90+0*(360/7)}:3.5) {};
            \node (345) at ({90+1*(360/7)}:3.5) {};
            \node (167) at ({90+2*(360/7)}:3.5) {};
            \node (234) at ({90+3*(360/7)}:3.5) {};
            \node (567) at ({90+4*(360/7)}:3.5) {};
            \node (123) at ({90+5*(360/7)}:3.5) {};
            \node (456) at ({90+6*(360/7)}:3.5) {};
            
            \node (245) at ({90+360/21+0*(360/7)}:4) {};
            \node (157) at ({90+360/21+1*(360/7)}:4) {};
            \node (134) at ({90+360/21+2*(360/7)}:4) {};
            \node (467) at ({90+360/21+3*(360/7)}:4) {};
            \node (237) at ({90+360/21+4*(360/7)}:4) {};
            \node (356) at ({90+360/21+5*(360/7)}:4) {};
            \node (126) at ({90+360/21+6*(360/7)}:4) {};
            
            \node (137) at ({90+2*360/21+0*(360/7)}:4) {};
            \node (346) at ({90+2*360/21+1*(360/7)}:4) {};
            \node (267) at ({90+2*360/21+2*(360/7)}:4) {};
            \node (235) at ({90+2*360/21+3*(360/7)}:4) {};
            \node (156) at ({90+2*360/21+4*(360/7)}:4) {};
            \node (124) at ({90+2*360/21+5*(360/7)}:4) {};
            \node (457) at ({90+2*360/21+6*(360/7)}:4) {};

            \draw (123) -- (456);
            \draw (123) -- (457);
            \draw (123) -- (467);
            \draw (123) -- (567);
            
            \draw (124) -- (356);
            \draw (124) -- (357);
            \draw (124) -- (567);
            
            \draw (125) -- (347);
            \draw (125) -- (367);
            
            \draw (126) -- (345);
            \draw (126) -- (457);
            \draw (126) -- (357);
            
            \draw (127) -- (456);
            \draw (127) -- (345);
            \draw (127) -- (346);
            \draw (127) -- (356);
            
            \draw (134) -- (257);
            \draw (134) -- (267);
            \draw (134) -- (567);
            
            \draw (135) -- (267);
            \draw (135) -- (467);
            
            \draw (136) -- (245);
            \draw (136) -- (457);
            
            \draw (137) -- (245);
            \draw (137) -- (246);
            \draw (137) -- (456);
            
            \draw (145) -- (236);
            \draw (145) -- (367);
            
            \draw (146) -- (235);
            \draw (146) -- (237);
            
            \draw (147) -- (236);
            \draw (147) -- (256);
            
            \draw (156) -- (234);
            \draw (156) -- (237);
            \draw (156) -- (247);
            
            \draw (157) -- (234);
            \draw (157) -- (246);
            \draw (157) -- (346);
            
            \draw (167) -- (234);
            \draw (167) -- (235);
            \draw (167) -- (245);
            \draw (167) -- (345);

            \draw (234) -- (567);
            \draw (235) -- (467);
            \draw (237) -- (456);
            \draw (247) -- (356);
            \draw (256) -- (347);
            \draw (257) -- (346);
            \draw (267) -- (345);

            \draw (136) to[bend left=25] (247);
            \draw (247) to[bend left=25] (135);
            \draw (135) to[bend left=25] (246);
            \draw (246) to[bend left=25] (357);
            \draw (357) to[bend left=25] (146);
            \draw (146) to[bend left=25] (257);
            \draw (257) to[bend left=25] (136);

            \draw (236) to[bend right =30] (457);
            \draw (256) to[bend right =30] (137);
            \draw (125) to[bend right =30] (346);
            \draw (145) to[bend right =30] (267);
            \draw (147) to[bend right =30] (235);
            \draw (347) to[bend right =30] (156);
            \draw (367) to[bend right =30] (124);
            
            \draw (236) to[bend left =30] (157);
            \draw (256) to[bend left =30] (134);
            \draw (125) to[bend left =30] (467);
            \draw (145) to[bend left =30] (237);
            \draw (147) to[bend left =30] (356);
            \draw (347) to[bend left =30] (126);
            \draw (367) to[bend left =30] (245);
            
            

            \foreach \x in {123,124,125,126,127,134,135,136,137,145,146,147,156,157,167,234,235,236,237,245,246,247,256,257,267,345,346,347,356,357,367,456,457,467,567}{

            \begin{scope}[shift={(\x)}]
                \filldraw[fill = white, thick] (0,0) circle (0.2);
                \foreach \y in {0,1,...,6}{
                \draw (0,0) -- ({90-360/7*\y}:0.2);
                }
                \filldraw[fill = gray!70] (0,0) -- ({90-360/7*(floor(\x/100)-1)}:0.2) arc({90-360/7*(floor(\x/100)-1)}:{90-360/7*floor(\x/100)}:0.2) -- cycle;
                \filldraw[fill = gray!70] (0,0) -- ({90-360/7*(floor(\x/10)-floor(\x/100)*10-1)}:0.2) arc({90-360/7*(floor(\x/10)-floor(\x/100)*10-1)}:{90-360/7*(floor(\x/10)-floor(\x/100)*10)}:0.2) -- cycle;
                \filldraw[fill = gray!70] (0,0) -- ({90-360/7*(\x-floor(\x/10)*10-1)}:0.2) arc({90-360/7*(\x-floor(\x/10)*10-1)}:{90-360/7*(\x-floor(\x/10)*10)}:0.2) -- cycle;
            \end{scope}}

            \foreach \x in {124,156,235,267,346,137,457,356,237,467,134,157,245,126}{

            \begin{scope}[shift={(\x)}]
                \filldraw[fill = white, thick] (0,0) circle (0.2);
                \foreach \y in {0,1,...,6}{
                \draw (0,0) -- ({90-360/7*\y}:0.2);
                }
                \filldraw[fill = red] (0,0) -- ({90-360/7*(floor(\x/100)-1)}:0.2) arc({90-360/7*(floor(\x/100)-1)}:{90-360/7*floor(\x/100)}:0.2) -- cycle;
                \filldraw[fill = red] (0,0) -- ({90-360/7*(floor(\x/10)-floor(\x/100)*10-1)}:0.2) arc({90-360/7*(floor(\x/10)-floor(\x/100)*10-1)}:{90-360/7*(floor(\x/10)-floor(\x/100)*10)}:0.2) -- cycle;
                \filldraw[fill = red] (0,0) -- ({90-360/7*(\x-floor(\x/10)*10-1)}:0.2) arc({90-360/7*(\x-floor(\x/10)*10-1)}:{90-360/7*(\x-floor(\x/10)*10)}:0.2) -- cycle;
            \end{scope}}
            
        \end{tikzpicture}
        \caption{A $\gamma_{\times 2}$-set of $\Kneser{7}{3}$}
        \label{fig: 2-upla en K(7,3)}
    \end{figure}
    It is enough to consider two disjoint Fano Planes $F_1,F_2$ and the sets 
    $$D_1=F_1,\qquad D_2=F_1\cup F_2,\qquad D_3=V\setminus D_2,\qquad D_4=V\setminus D_1,\qquad D_5=V,$$
    where $V=V(\Kneser{7}{3})$. For each $k\in[5]$, $D_k$ turns out to be a $k$-tuple dominating set with cardinality $|D_k|=7k$ and in consequence is a   $\gamma_{\times k}$-set.

\end{remark}

In Remark \ref{remarkkdom73} we provide $\gamma_{\times k}$-sets of $\Kneser{7}{3}$. In fact, by exhaustive methods we can prove that these sets are unique up to automorphism. Regarding Steiner systems $S(4,5,11)$, it is also known that there are exactly two disjoint such Steiner systems. We know that the union of them is a $\gamma_{\times 2}$-set of $\Kneser{11}{5}$. We also find, for $k\in\{3,4,5,6\}$, a $\gamma_{\times k}$-set of $\Kneser{11}{5}$ (see Appendix \ref{sec: K(11,5)}). Although this parallelism with the case $\Kneser{7}{3}$, we do not known if the given $\gamma_{\times k}$-sets of $\Kneser{11}{5}$ are unique up to automorphism for $k\in\{2,3,4,5\}$.

\section{\texorpdfstring{General results on $k$-tuple domination on Kneser graphs}{General results on k-tuple domination on Kneser graphs}}\label{sec:general}

We note that the bounds for the $k$-tuple domination number given in Theorem \ref{thm: cota packing} are not tight in general. For example, as we have mentioned, $\rho(\Kneser{n}{r})=1$ if $n\geq3r-1$, but $\kdom{n,r}{1}=\dom{\Kneser{n}{r}}\geq r+1$ for $n\geq r^2+r$. Also, the known general bounds for $k$-tuple domination number of graphs \cite{martinez2022note,chang2008upper,gagarin2008generalised,RAUTENBACH200798} are not efficient to determine the $\kdom{n,r}{k}$ for $n\geq2r+1$. Moreover, some of these upper bounds applied to $\kdom{n,r}{k}$ are increasing with respect to $n$ whereas, as we will prove in Theorem \ref{thm: monotonia respecto a n}, the parameter $\kdom{n,r}{k}$ is decreasing with respect to $n$. Therefore, in this section we will apply properties of Kneser graphs to study $k$-tuple dominating sets on $\Kneser{n}{r}$ for $n$ large enough. 

Since Kneser graphs are regular graphs and each vertex in $\Kneser{n}{r}$ has degree $\binom{n-r}{r}$, then the Kneser graph $\Kneser{n}{r}$ can only have a $k$-tuple dominating set if $k\leq \binom{n-r}{r}+1$. Furthermore, if $k=\binom{n-r}{r}+1$, the entire set of vertices in $\Kneser{n}{r}$ is the only $k$-tuple dominating set. 
The Kneser graph $\Kneser{n}{1}$ is isomorphic to the complete graph $K_n$. In this case, $\kdom{n,1}{k}=k$ for $k\in[n]$.
On the other hand, if $n=2r$, then the Kneser graph $\Kneser{2r}{r}$ is isomorphic to $\frac{1}{2}\binom{2r}{r}$ copies of $K_2$. It is noteworthy that $\kdom{2r,r}{1}=\frac{1}{2}\binom{2r}{r}$, $\kdom{2r,r}{2}=\binom{2r}{r}$, and no $k$-tuple dominating set exists for these graphs if $k\geq 3$. 
Thus, for the subsequent discussion, we consider $n\geq 2r+1$ and $r\geq 2$.

{In Theorem \ref{TheoDomBigN} it is shown that if $n\geq r^2+r$, then $\kdom{n,r}{1}=\dom{\Kneser{n}{r}}{}=r+1$. Moreover, for every $k$ we state the following result.}

\begin{lemma}
\label{lem: k+r para n grandes}
If $n\geq r(k+r)$, then $\kdom{n,r}{k}=k+r$.
\begin{proof}
    Let $D=\{u_1,\ldots,u_{k+r}\}$ be a set of vertices of $\Kneser{n}{r}$ such that $u_i\cap u_j=\emptyset$ for all $i\neq j$. It is possible since $n\geq r(k+r)$. Consider a vertex $u$ of $\Kneser{n}{r}$. If $u\in D$, then there exists $i\in\{1,\ldots,k+r\}$ such that $u=u_i$. Since $u\cap u_j=\emptyset$ for all $j\neq i$, we have $D\setminus\{u_i\}\subseteq N(u)$. Thus, 
    $$|D\cap N[u]|=|D|=k+r\geq k.$$
    If $u\notin D$, we have $|\{i:u_i\cap u\neq\emptyset\}|\leq r$. So
    $$|D\cap N[u]|\geq|D|-r=k.$$    
    Therefore $D$ is a $k$-tuple dominating set of $\Kneser{n}{r}$ with cardinality $k+r$, and  $\kdom{n,r}{k}\leq k+r$.
        
    Now, assume that there exists a $k$-tuple dominating set $D$ with cardinality $|D|=k+r-1$. Let us consider $r$ distinct vertices $u_1,\ldots,u_r$ of $D$. Let $a_1\in u_1$ and for each $2\leq i\leq r$ we choose $a_i\in u_i$ such that $a_i\notin\{a_1,\ldots,a_{i-1}\}$. Let $w=\{a_1,\ldots,a_r\}$. If $w\neq u_i$ for every $i\in [r]$, then we have $w\cap u_i\neq \emptyset$, and so $u_i\notin N[w]$ for all $i\in [r]$. It turns out that
    $$|D\cap N[w]|\leq|D|-r\leq k-1,$$
    contradicting the fact that $D$ is $k$-tuple dominating.
    
    If $w=u_j$ for some $j$, note that we can choose $b\in [n]\setminus\left(\bigcup_{i=1}^{r}u_i\right)$. This choice is possible since $|\bigcup_{i=1}^r u_i|\leq n-1$. In fact, $u_j\subseteq \left(\bigcup_{i\neq j}u_i\right)\cup\{a_j\}$ and then
    $$\left|\bigcup_{i=1}^r u_i \right|\leq \left|\bigcup_{i\neq j} u_i \right|+1 \leq r(r-1)+1<r(k+r)\leq n.$$
    Let $w'=\left[w\setminus\{a_j\}\right]\cup\{b\}$. We have $w'\neq u_i$ for every $i\in [r]$. Besides, $a_i\in w'\cap u_i$ for $i\neq j$, and $w'\cap u_j=w\setminus\{a_j\}$. Therefore $w'\cap u_j\neq \emptyset$. Thus, $u_i\notin N[w']$ for all $i\in [r]$. Then
    $$|D\cap N[w']|\leq|D|-r\leq k-1,$$
    contradicting the fact that $D$ is $k$-tuple dominating. 
    
    Therefore, it cannot exist a $k$-tuple dominating set of $\Kneser{n}{r}$ with cardinality less than $k+r$. We conclude that $\kdom{n,r}{k}=k+r$.
\end{proof}
\end{lemma}

Note that the condition $n\geq r(k+r)$ in Lemma \ref{lem: k+r para n grandes} guarantees the existence of a set of $k+r$ pairwise disjoint vertices of $\Kneser{n}{r}$.
In the next result, we state that unless $k=1$ and $r=2$, this kind of set families characterizes the $\gamma_{\times k}$-sets in $\Kneser{n}{r}$ for these values of $n$.

\begin{lemma}
\label{lem: k+r vertices disjuntos}
Let $k\geq 2$ if $r=2$ and $k\geq 1$, otherwise. If $D$ is a $k$-tuple dominating set of $\Kneser{n}{r}$ with cardinality $k+r$, then the vertices of $D$ are pairwise disjoint.

\begin{proof}
Let $D$ be a $k$-tuple dominating set of $\Kneser{n}{r}$ with cardinality $k+r$. Assume that the vertices in $D$ are not pairwise disjoint. Let $a\in [n]$ such that $i_a=\max_{x\in[n]}i_x$. Under our assumption, $i_a\geq 2$. 

If $i_a\geq 3$, let $u_1,u_2,u_3$ in $D$ such that $a\in u_1\cap u_2\cap u_3$. Let us consider $u_4,\ldots, u_{r+1}$ vertices in $D\setminus\{u_1,u_2,u_3\}$. Let $b_1=a$. For $j=2,\ldots,r-1$, we choose $b_j\in u_{j+2}\setminus\{b_1,\ldots,b_{j-1}\}$. Let $b=\{b_1,\ldots,b_{r-1}\}$. Since $n\geq 2r+1$, it follows that $|[n]\setminus b|=n-(r-1)\geq r+2$. Thus, there exists $x\in [n]$ such that $w=b\cup \{x\}\neq u_j$ for $j=1,\ldots,r+1$. As a result, for every $j$ we have $w\cap u_j\neq \emptyset$ and $w\neq u_j$. Then,
\begin{equation}
\label{eq: k+r disj 1}
    |D\cap N[w]|\leq |D|-|\{u_1,\ldots,u_{r+1}\}|=|D|-(r+1)=k-1,
\end{equation}
which contradicts the fact that $D$ is a $k$-tuple dominating set. So, $i_a=2$.

Let $u_1$ and $u_2$ be the two vertices that contain the element $a$. 

If $r>2$, let $u_3,\ldots,u_{r+1}$ be vertices in $D\setminus\{u_1,u_2\}$. Let $b_1=a$ and for $j=2,\ldots,r$, we choose $b_j\in u_{j+1}\setminus\{b_1,\ldots,b_{j-1}\}$. Let $w=\{b_1,\ldots,b_{r}\}$. $b$ is a vertex that is not adjacent to any $u_j$ for $j=1,\ldots,r+1$. If $w\neq u_j$ for each $j$, then it holds (\ref{eq: k+r disj 1}), contradicting the fact that $D$ is $k$-tuple dominating. Then $w=u_j$ for some $j$. Since the only vertices that contain the element $b_1=a$ are $u_1$ and $u_2$, without loss of generality, $w=u_1$. Then we have $b_j\in u_1\cap u_{j+1}$ for every $j=1,\ldots,r$. Let $x\in u_2\setminus w$ and let $w'=w\setminus\{b_1\}\cup\{x\}=\{x,b_2,\ldots,b_{r}\}$. Let us see that $w'\neq u_j$ and $w'\cap u_j\neq \emptyset$ for every $j$. In fact, as $b_j\in u_1\cap u_{j+1}$ for $j=2,\ldots,r$, then $w'\cap u_j\neq \emptyset$ for every $j\in 1,\ldots,r+1$. Moreover, if $w'=u_\ell$ for some $\ell$, then since $a\notin w'$, $\ell\geq 3$ and for $j\neq \ell$, $j\geq 2$ we have that $b_{j-1}\in u_1\cap u_j\cap u_\ell$. Thus, $i_{b_{j-1}}\geq 3$ but $\max_{x\in[n]}i_x=2$. Therefore, $w'\neq u_j$ for every $j$, and in consequence it holds (\ref{eq: k+r disj 1}) for $w'$.

If $r=2$, we have $u_1=\{a,c\}$, $u_2=\{a,d\}$ for some $c,d\in [n]$. Since $k\geq 2$, there exists a vertex $u_3\in D\setminus\{u_1,u_2\}$ such that $u_3\neq\{c,d\}$. Let $x\in u_3\setminus\{c,d\}$. Since $i_a=2$, $x\neq a$. Consider the vertex $w=\{a,x\}$. We have $w\neq u_j$ and $w\cap u_j\neq\emptyset$ for every $j$. Then it holds (\ref{eq: k+r disj 1}).

Either if $r>2$ or $r=2$ we arise to a contradiction since $D$ is a $k$-tuple dominating set.

Therefore, we conclude that the vertices in $D$ are pairwise disjoint.
\end{proof}
\end{lemma}

This result does not hold when $k=1$ and $r=2$. In fact, in \cite{ivanvco1993domination} it is shown that the dominating sets of $\Kneser{n}{2}$ for $n\geq 5$ are the sets of $3$ vertices that are either pairwise disjoint or mutually intersecting. 

As a by product of lemmas \ref{lem: k+r para n grandes} and \ref{lem: k+r vertices disjuntos} we have the following result for $n$ large enough.

\begin{theorem}
\label{thm: kdom = k+r}
    For $k\geq 2$, $\kdom{n,r}{k}=k+r$ if and only if $n\geq r(k+r)$.
\end{theorem}

For $n<r(k+r)$ we will see that $\kdom{n,r}{k}$ is greater than $k+r$. To this end we observe monotonicity properties of the parameter $\kdom{n,r}{k}$ with respect to $k$ and $n$.

Note that if $D$ is a $(k+1)$-tuple dominating set of a graph $G$ and $v\in D$, then $D\setminus\{v\}$ is a $k$-tuple dominating set ($k\geq1$).
Therefore, for any graph $G$, we have $\kdom{G}{k}<\kdom{G}{(k+1)}$ whenever $\mindeg{G}\geq k-1$. Then it turns out that for fixed $r$ and $n$, and $2 \leq k\leq \binom{n-r}{r}$, we have $\kdom{n,r}{k} < \kdom{n,r}{(k+1)}$, i.e. $\kdom{n,r}{k}$ is strictly increasing with respect to $k$. 

In addition, notice that there is not monotonicity for the $k$-tuple domination number with respect to induced subgraphs. However, in the case of Kneser graphs, if $2r\leq m<n$, $\Kneser{m}{r}$ is an induced subgraph of $\Kneser{n}{r}$ and we will show that $\kdom{m,r}{k}\geq \kdom{n,r}{k}$ for every $k\geq 2$ such that $k\leq\binom{m-r}{r}+1$. This is, for each fixed $k\geq 2$, $\kdom{n,r}{k}$ is decreasing with respect to $n$. 

\begin{theorem}
\label{thm: monotonia respecto a n}
If $D$ is a $k$-tuple dominating set of $\Kneser{n}{r}$ with $k\geq 2$, then $D$ is a $k$-tuple dominating set of $\Kneser{n+1}{r}$. In consequence, we have $\kdom{n+1,r}{k}\leq \kdom{n,r}{k}$ for every $n$.

\begin{proof} 
Let $D$ be a $k$-tuple dominating set of $\Kneser{n}{r}$. In order to see that $D$ is a $k$-tuple dominating set of $\Kneser{n+1}{r}$, it is enough to show that every vertex $u$ of $\Kneser{n+1}{r}$ containing the element $n+1$ satisfies
$$|N_{n+1}[u]\cap D|\geq k,$$
where $N_{t}[u]$ denotes the neighbourhood of the vertex $u$ in the graph $\Kneser{t}{r}$.

Let $\tilde{u}$ be a $(r-1)$-subset of $[n]$ and $u=\tilde{u}\cup\{n+1\}$. Let us define, for $b\in [n]\setminus \tilde{u}$, $u_b=\tilde{u}\cup \{b\}$. Note that if $b\in [n]\setminus \tilde{u}$, then
$$|N_{n+1}[u]\cap D|=|\{w\in D:|w\cap \tilde{u}|=0\}|\geq  |N_{n}(u_b)\cap D|\geq \begin{cases}
    k-1, &\text{if }u_b\in D\\
    k, &\text{if }u_b\notin D
\end{cases}.$$
Therefore, if there exists $b\in [n]\setminus \tilde{u}$ such that $u_b\notin D$, then $|N_{n+1}[u]\cap D|\geq k$.

Otherwise, consider $u_b$ for some $b\in [n]\setminus\tilde{u}$. We have
$u_b\in D$ and in consequence, $|N_{n}(u_b)\cap D|\geq k-1$. Consider $z\in N_{n}(u_b)\cap D$, and let $x$ be an element from $z$. Now, consider the vertex $u_x$. Note that $w\cap \tilde{u}=\emptyset$ for every $w\in N_n(u_x)\cap D$. Besides, $\tilde{u}\cap z=\emptyset$ since $z\in N_n(u_b)$. It follows
$$|N_{n+1}[u]\cap D|=|\{w\in D:|w\cap \tilde{u}|=0\}|\geq  |N_{n}(u_x)\cap D|+1\geq k.$$

Therefore $D$ is a $k$-tuple dominating set of $\Kneser{n+1}{r}$.

\end{proof}
\end{theorem}

From theorems \ref{thm: kdom = k+r} and \ref{thm: monotonia respecto a n} we have the following.

\begin{corollary}
    For $k\geq 2$, $\kdom{n,r}{k}\geq k+r$. Moreover, if $n<r(k+r)$, then $\kdom{n,r}{k}\geq k+r+1$.
\end{corollary}

The remaining of this section is devoted to obtain $\gamma_{\times k}$-sets for $k=\binom{n-r}{r}-t$ when $n$ is large enough with respect to both $r$ and $t$.
As we have mentioned, if $k=\binom{n-r}{r}+1$, the only $\gamma_{\times k}$-set of $\Kneser{n}{r}$ is the set of vertices itself. On the other hand, when $n\geq 3r-1$, the diameter of the Kneser graph $\Kneser{n}{r}$ is equal to $2$ \cite{valencia2005diameter}. In these cases, any pair of vertices in $\Kneser{n}{r}$ are adjacent or they have a common neighbour. Thus, for $k=\binom{n-r}{r}$ and $n\geq 3r-1$, it follows that $\kdom{n,r}{k}=\binom{n}{r}-1$.
In a similar way, we prove that for a positive integer $t$, when $n\geq (t+3)r-\left\lceil\frac{t+2}{2}\right\rceil$, every set of $t+2$ vertices in $\Kneser{n}{r}$ is contained in the closed neighborhood of some vertex and in consequence, for $k=\binom{n-r}{r}-t$, it follows that $\kdom{n,r}{k}=\binom{n}{r}-(t+1)$.

\begin{remark}\label{rem: cardinal de la union}
    Let $t\in\mathbb{N}$, $t\geq 2$, and let $S$ be a set of $t$ vertices of $\Kneser{n}{r}$ such that each vertex $u\in S$ intersects at least one vertex in $S\setminus\{u\}$. It holds
    \begin{equation} \label{HI}
            \left|\bigcup_{v\in S} v\right|\leq tr-\left\lceil\frac{t}{2}\right\rceil
        \end{equation}
\end{remark}

\begin{lemma}\label{lem: t vertices}
    Let $t\geq 2$. If $n\geq (t+1)r-\left\lceil\frac{t}{2}\right\rceil$ and $S$ is a set of $t$ vertices of $\Kneser{n}{r}$, then there exists a vertex $w$ of $\Kneser{n}{r}$ such that $S\subseteq N[w]$.
    \begin{proof}
        Let $t,r,n\in \mathbb{N}$ such that $t\geq 2$ and $n\geq (t+1)r-\left\lceil\frac{t}{2}\right\rceil$, and let $S$ be a set of $t$ vertices of $\Kneser{n}{r}$. If there exists $w\in S$ such that $w\cap u=\emptyset$ for every $u\in S\setminus\{w\}$, the result holds. Now, assume that each vertex of $S$ intersects at least another vertex in $S$. Thus, by Remark \ref{rem: cardinal de la union}, it follows (\ref{HI}).
        Since $n\geq (t+1)r-\left\lceil\frac{t}{2}\right\rceil$, we have $\left|[n]\setminus \bigcup_{v\in S} v\right|=n-\left| \bigcup_{v\in S} v\right|\geq r$ and there exists at least one vertex $w\in [n]\setminus\bigcup_{v\in S}v$. Since $w\cap v=\emptyset$ for each $v\in S$, we have $S\subseteq N[w]$, and the statement holds.
    \end{proof}
\end{lemma}

\begin{theorem}\label{thm: large k}
    For $t\in \mathbb{N}\cup\{0\}$, $n\geq (t+3)r-\left\lceil\frac{t+2}{2}\right\rceil$ and $k=\binom{n-r}{r}-t$, it holds $\kdom{n,r}{k}=\binom{n}{r}-(t+1)$.
    Moreover, for every set $S\subseteq V(\Kneser{n}{r})$ with cardinality $t+1$, $V(\Kneser{n}{r})\setminus S$ is a $\gamma_{\times k}$-set.

    \begin{proof}
        Let $S$ be a set of vertices of $\Kneser{n}{r}$ with cardinality $t+2$. Let us show that the set $D=V(\Kneser{n}{r})\setminus S$ is not a $k$-tuple dominating set. In fact, by Lemma \ref{lem: t vertices}, we have that there exists a vertex $w\in V(\Kneser{n}{r})$ such that $S\subseteq N[w]$. Thus, $|N[w]\cap D|=|N[w]|-|S|=\binom{n-r}{r}+1-(t+2)=k-1$. In consequence, $V(\Kneser{n}{r})\setminus S$ is not a $k$-tuple dominating set of $\Kneser{n}{r}$. Since it is true for any set of vertices $S$ with cardinality $t+2$, it follows that $\kdom{n,r}{k}\geq \binom{n}{r}-(t+1)$.

        On the other hand, let $S$ be any set of vertices of $\Kneser{n}{r}$ with cardinality $t+1$, and $D=V(\Kneser{n}{r})\setminus S$. Let $u\in V(\Kneser{n}{r})$. We have
        $$|N[u]\cap D|=|N[u]|-\underbrace{|N[u]\cap S|}_{\leq |S|}\geq \binom{n-r}{r}+1-(t+1)=k$$
        Thus, $D$ is a $k$-tuple dominating set, and $\kdom{n,r}{k}\leq |D|=\binom{n}{r}-(t+1)$.

        Therefore, $\kdom{n,r}{k}=\binom{n}{r}-(t+1)$.
    \end{proof}
\end{theorem}

\begin{remark}
    Note that if $n=(t+3)r-\left\lceil\frac{t+2}{2}\right\rceil-1$ and $k=\binom{n-r}{r}-t$, then $\kdom{n,r}{k}<\binom{n}{r}-(t+1)$. In fact, we may consider the following set $S$. For $t$ even,
    \begin{align*}
        S = \bigg\{&[\xi+1..\xi+r],[\xi+r..\xi+2r-1],\,\xi= (x-1)(2r-1), \text{ with }x\in \left[\frac{t}{2}+1\right] \bigg\}
    \end{align*}
    and for $t$ odd,
    \begin{align*}
        S = \bigg\{&[\xi+1..\xi+r],[\xi+r..\xi+2r-1],\,\xi= (x-1)(2r-1), \text{ with }x\in \left[\left\lceil\frac{t}{2}\right\rceil+1\right] \bigg\}\cup\\
        \bigg\{&[\xi+1..\xi+r],[\xi+r..\xi+2r-1],[\xi+2r..\xi+3r-2],\,\xi= \left\lceil\frac{t}{2}\right\rceil(2r-1)\bigg\}.
    \end{align*}
    In both cases, $S$ is a set of $t+2$ vertices and $V(\Kneser{n}{r})\setminus S$ is a $k$-tuple dominating set. Thus, $\kdom{n,r}{k}\leq \binom{n}{r}-(t+2)$.
\end{remark}

\section{\texorpdfstring{$k$-tuple domination on Kneser graphs $\Kneser{n}{2}$}{k-tuple domination on Kneser graphs Kn(n,2)}}\label{secKn2}

The study of subsets of vertices satisfying certain restriction in Kneser graphs $\Kneser{n}{2}$ deserves a particular analysis, since this subclass of Kneser graphs has a remarkable structure. For instance, in \cite{bresar2019grundy,qAnalog,GenPosKneser,TreewidthKneser,jalilolghadr2023total,ivanvco1993domination} the authors studied the corresponding invariant in the subclass $\Kneser{n}{2}$ separately. As we have mentioned, in \cite{ivanvco1993domination} it is shown that $\kdom{n,2}{1}=3$ for every $n\geq4$. In this section we focus on the $k$-tuple domination number $\kdom{n,2}{k}$ for $k\geq 2$.
To this end, let us see a characterization of the $k$-tuple dominating sets of $\Kneser{n}{2}$ in terms of the occurrences of the elements in $[n]$.

\begin{lemma}
    \label{lemma: cotas indices}
    Let $n\geq 5$. $D$ is a $k$-tuple dominating set of $\Kneser{n}{2}$ if and only if for every pair $a,b\in[n]$
    $$i_a+i_b\leq
    \begin{cases}
    |D|-k+2, & \text{if }\{a,b\}\in D,\\
    |D|-k, & \text{if }\{a,b\}\notin D.
    \end{cases}$$
    
    \begin{proof}
        Let $D$ be a set of vertices of $\Kneser{n}{2}$ and let $u=\{a,b\}$ be a vertex of $\Kneser{n}{2}$. The amount of vertices in $D$ that contain either the element $a$ or $b$ is exactly $i_a+i_b-1$ when $u\in D$ and $i_a+i_b$ otherwise. 
        
        On the one hand, if $u\in D$, then we have
        \begin{equation}
        \label{eq: cotas indices 1}
            \begin{aligned}
                |D\cap \neighc{u}{}|&=|D|-|\{v\in D:v\neq u \wedge v\cap u\neq\emptyset\}|=\\
                &=|D|-\left(|\{v\in D: a\in v\vee b\in v\}|-1\right)=\\
                &=|D|-(i_a+i_b-1)+1= |D|+2-(i_a+i_b).
            \end{aligned}
        \end{equation}
        
        On the other hand, whether $u\notin D$ it holds that
        \begin{equation}
        \label{eq: cotas indices 2}
            |D\cap \neighc{u}{}|=|D|-|\{v\in D: a\in v\vee b\in v\}|=|D|-(i_a+i_b)= |D|-(i_a+i_b).
        \end{equation}
        
        From (\ref{eq: cotas indices 1}) and (\ref{eq: cotas indices 2}) it turns out that $D$ is a $k$-tuple dominating set if and only if for every pair of elements $a,b\in [n]$ it holds
        $$\begin{cases}
        |D|+2-(i_a+i_b)\geq k, & \text{if }\{a,b\}\in D,\\
        |D|-(i_a+i_b)\geq k, & \text{if }\{a,b\}\notin D.
        \end{cases}$$
        Or, equivalently
        $$i_a+i_b\leq
        \begin{cases}
        |D|-k+2, & \text{if }\{a,b\}\in D,\\
        |D|-k, & \text{if }\{a,b\}\notin D.
        \end{cases}$$
    \end{proof}
\end{lemma}

\begin{remark}\label{rem: indices}
    Note that if $D$ is a $k$-tuple dominating set and $\{a,b\}\in D$, then from Lemma \ref{lemma: cotas indices} together with the fact that $i_b\geq 1$, we have that for any $a\in [n]$, it holds $i_a\leq |D|-k+1$.
\end{remark}

From Theorem \ref{thm: kdom = k+r} with $r=2$, we have that $\kdom{n,2}{k}=k+2$ if and only if $n\geq 2(k+2)$. In addition, from monotonicity, it follows that for $k\geq 2$ and $n< 2(k+2)$, $\kdom{n,2}{k}\geq k+3$. Moreover, in the following result we state that the only value of $n$ for which this bound is tight is $n=2k+3$. 
We introduce the following notation. Given a set $A$ and a positive integer $r\leq |A|$, we denote by $\binom{A}{r}$ the set of all the $r$-subsets of $A$.

\begin{theorem}\label{thm: kdom(n,2)=k+3}
    If $k\geq 2$, then $\kdom{n,2}{k}=k+3$ if and only if $n=2k+3$.

    \begin{proof}

        Let $k\geq 2$ and $n\leq 2k+3$. Let us suppose that there exists a $k$-tuple dominating set $D=\{u_1,\ldots,u_{k+3}\}$ of $\Kneser{n}{2}$ with cardinality $k+3$. We will prove that $n$ must be equal to $2k+3$. Note that it is enough to prove that $n\geq 2k+3$.
        Let us consider a vertex $u$, and let $a\in u$. By Remark \ref{rem: indices} we have $i_a\leq 4$. Moreover, we can prove the following.

\medskip

\begin{claim}
    $i_a\leq 2$ for every $a\in [n]$.
\end{claim}

\begin{proofclaim}
    Suppose, on the contrary, that $i_a\geq 3$ for some $a\in[n]$. 
    
    Note that there exists $b\in [n]-\{a\}$ such that $\{a,b\}\notin D$ and $i_b\geq 1$.
    In fact, if for every element $x\neq a$ with $i_x\geq 1$ we have $\{a,x\}\in D$, then there are at most $i_a+1$ elements in $[n]$ that appear in vertices of $D$. Let $x$ be an element different from $a$ such that $i_x\geq 1$. It turns out that $\{a,x\}\in D$ and by Lemma \ref{lemma: cotas indices} $i_a+i_x\leq 5$.
    
    Then if $i_a=4$, we have $i_x=1$ and $k+3=|D|=4$ which contradicts that $k\geq 2$. On the other hand, if $i_a=3$ then it turns out that $i_y\leq 2$ for every $y\neq a$ with $i_y\geq 1$, and we have
    $$2(k+3)=2|D|=\sum_{y=1}^{n}i_y\leq i_a+2i_a=9.$$
    We also get $k+3=4$, which does not hold since $k\geq 2$.
    
    Therefore, it is possible to choose $b\in[n]$ such that $i_b\geq 1$ and $\{a,b\}\notin D$. Thus $i_a\geq 3$, $i_b\geq 1$, and in consequence $i_a+i_b\geq 4$, which leads to a contradiction since by Lemma \ref{lemma: cotas indices} we have $i_a+i_b\leq 3$. 
    
    We conclude $i_a\leq 2$ for all $a\in [n]$.
\end{proofclaim}

Let $n_1=|X_1|$ and $n_2=|X_2|$. We have $n\geq n_2+n_1$ and $2|D|=2n_2+n_1$. So
$$n\geq n_2+n_1=2|D|-n_2=2(k+3)-n_2=2(k+2)+(2-n_2).$$
Since $n<2(k+2)$, it turns out that $n_2\geq 3$. Let us see that $n_2$ is exactly $3$.

\medskip

\begin{claim}
    $n_2=3$.
\end{claim}

\begin{proofclaim}
    Suppose that $a$ and $b$ are two elements in $[n]$ such that $i_a=i_b=2$, and let $v=\{a,b\}$. If $v\notin D$, by Lemma \ref{lemma: cotas indices} $i_a+i_b\leq 3$ which leads to a contradiction. Thus, $v\in D$. 
    If $n_2\geq 4$, let us consider $a_1,a_2,a_3,a_4$ such that $i_{a_j}=2$ for $1\leq j \leq 4$. Then $\{a_1,a_j\}\in D$ for every $2\leq j\leq 4$, and $i_{a_1}\geq 3$ which contradicts the fact that $i_{a_1}=2$. Thus, $n_2=3$.
\end{proofclaim}

Therefore we have $n\geq n_2+n_1=2k+3$. In consequence, it turns out that if $\kdom{n,2}{k}=k+3$, then $n=2k+3$.

On the other hand, if $n=2k+3$, then from Theorem \ref{thm: kdom = k+r} and monotonicity, we have $\kdom{n,2}{k}\geq k+3$. Now, let us see that there exists a $k$-tuple dominating set $\widehat{D}$ with cardinality $k+3$. From the reasoning above, we have that $\widehat{D}=X_1\cup X_2$, $|X_2|=3$ and $\binom{X_2}{2}\subseteq\widehat{D}$. 
In fact, consider the set $\widehat{D}$ given by
$$\widehat{D}=\binom{[3]}{2}\cup \{\{2a,2a+1\}:2\leq a\leq k+1 \}.$$
$\widehat{D}$ is a $k$-tuple dominating set of $\Kneser{n}{2}$ with $|\widehat{D}|=k+3$. Therefore $\kdom{2k+3,2}{k}=k+3$.

We conclude that for $k\geq 2$, $\kdom{n,2}{k}=k+3$ if and only if $n=2k+3$.

\end{proof}

\end{theorem}

As a by-product of the proof of Theorem $\ref{thm: kdom(n,2)=k+3}$ we have that if $n=2k+3$ with $k\geq 2$, and $D$ is a $\gamma_{\times k}$-set of $\Kneser{n}{2}$, 
then up to automorphism
$$D= \binom{[3]}{2} \cup \{\{2a,2a+1\}:2\leq a\leq k+1 \}.$$
Notice that in theorems \ref{thm: kdom = k+r} and \ref{thm: kdom(n,2)=k+3}, we provide $\kdom{n,2}{k}$ for $n\geq 2k+3$. From monotonicity on $n$ it follows that if $n\leq 2k+2$, then $\kdom{n,2}{k}\geq k+4$. Moreover, in the remaining of this section we will obtain $\kdom{n,2}{k}$ if $n$ is large enough with respect to $k$. To this end we give lower bounds for $\kdom{n,2}{k}$, and afterwards we prove that these bounds are tight by providing $k$-tuple dominating sets of $\Kneser{n}{2}$ with the desired cardinality.

\begin{proposition}\label{prop: lower bounds kdom(n,2) V2}
    Let $\alpha,n\in\mathbb{N}$ with $\alpha\geq 2$ and $n\geq 2\alpha+3+(\alpha \mod 2)$
    \begin{enumerate}
        \item $\kdom{n,2}{k}\geq k+2\alpha$, if $\frac{2}{\alpha}k+4\leq n<\frac{2}{\alpha-1}k+3$;
        \item $\kdom{n,2}{k}\geq k+2\alpha+1$, if $n=\left\lceil\frac{2}{\alpha}k\right\rceil+3$.
    \end{enumerate}

    \begin{proof}
        Let $\alpha,n\in \mathbb{N}$ with $\alpha\geq 2$ and $n\geq 2\alpha+3+(\alpha \mod 2)$. First, we prove the two following claims.

        \begin{claim}
            If $\Kneser{n}{2}$ admits a $k$-tuple dominating set with cardinality $k+2\alpha$, then $\alpha n\geq 2k+4\alpha$.
        \end{claim}
        
        \begin{proofclaim}
            Let $D$ be a $k$-tuple dominating set of $\Kneser{n}{2}$ with cardinality $|D|=k+2\alpha$. Our goal is to prove that $\alpha n\geq 2(k+2\alpha)=2|D|$. Let us suppose, on the contrary, that $\alpha n < 2|D|$.
            
            On the one side, by Lemma \ref{lemma: cotas indices} for every pair of elements $a,b\in[n]$ we have
            \begin{equation}
            \label{eq: kupla general rayo grande}
                i_a+i_b\leq\begin{cases}
                |D|-k+2=2(\alpha+1) & \text{if }\{a,b\}\in D,\\
                |D|-k=2\alpha & \text{if }\{a,b\}\notin D.
                \end{cases}
            \end{equation}
            Since $\alpha n < 2|D|=\sum_{x\in [n]}i_x$, we have that there is at least one element $a\in [n]$ for which $i_a\geq \alpha+1$. For every other element $b\in X_\alpha^\geq$, it holds $i_a+i_b\geq 2\alpha +1$, and by (\ref{eq: kupla general rayo grande}) $\{a,b\}\in D$.
            
            If $|X_\alpha^\geq|>1$, then we have that for every $x\in X_\alpha^\geq$, $i_x\leq \alpha +2$. In fact, if for some element $y$ it holds $i_y>\alpha +2$, then by (\ref{eq: kupla general rayo grande}) for every other $x\in [n]$
            $$i_x\leq 2(\alpha +1)-i_y<2(\alpha+1)-(\alpha+2)=\alpha,$$
            and $y$ would be the only element in $X_\alpha^\geq$, but $|X_\alpha^\geq|>1$. 
            
            Let us see that $|X_\alpha^\geq|>1$. Otherwise, $X_\alpha^\geq=\{a\}$ with $i_a\geq \alpha+1$ and in consequence
            \begin{align*}
                \alpha n < 2|D| &= \sum_{x\in [n]} i_x = i_a + \sum_{x\in [n]\setminus\{a\}} \overbrace{i_x}^{\leq \alpha-1} \leq i_a + (n-1)(\alpha - 1) = i_a + \alpha n - n - \alpha + 1.
            \end{align*}
            Then, $i_a > n + (\alpha - 1) \geq n$, which cannot be true since $i_x\leq n-1$ for every element $x\in [n]$. Therefore, $|X_\alpha^\geq|>1$. 
            In consequence, $i_x\leq \alpha +2$ for every $x\in X_\alpha^\geq$. What is more, if $y\in X_{\alpha +2}$, then for every $x\in X_{\alpha}^\geq$ 
            we have $i_x\leq 2(\alpha +1)-i_y=\alpha$ and it turns out that $i_x=\alpha$. Then either $X_\alpha^\geq=X_\alpha\cup X_{\alpha+1}$ or $X_\alpha^\geq=X_\alpha\cup X_{\alpha+2}$ with $|X_{\alpha+2}|=1$.
            
            \begin{description}
                \item[Case 1.] $X_\alpha^\geq =X_\alpha\cup X_{\alpha+1}$.
                
                Since for $a\in X_{\alpha+1}$ and every $b\in X_\alpha^\geq$ we have $\{a,b\}\in D$, then $|X_\alpha^\geq|\leq i_a+1 = \alpha+2$, and we have
                \begin{align*}
                    2|D| &= \sum_{x\in [n]} i_x = \sum_{x\in X_\alpha^\geq} i_x + \sum_{x\in X_{\alpha-1}^\leq} i_x\leq\\
                    &\leq |X_{\alpha}^\geq|(\alpha +1)+\underbrace{|X_{\alpha-1}^\leq|}_{n-|X_\alpha^\geq|}(\alpha -1)=\\
                    &=n(\alpha-1)+2|X_\alpha^\geq|\leq n\alpha -n + 2(\alpha+2).
                \end{align*}
                Since $2|D|\geq\alpha n+1$, it turns out that $n\leq 2(\alpha+2)-1$. By hypothesis, we have
                $$n\geq \begin{cases}
                    2(\alpha+2), &\text{if $\alpha$ is odd},\\
                    2(\alpha+2)-1, &\text{if $\alpha$ is even}.
                \end{cases}$$
                Thus, $\alpha$ is even and $n=2(\alpha+2)-1$. Therefore
                \begin{equation*}
                    \alpha n +1\leq 2|D|\leq n\alpha -n + 2(\alpha+2)=\alpha n +1.
                \end{equation*}
                It turns out that $2|D|=\alpha n +1$. But $2|D|$ is even whereas $\alpha n +1$ is odd.
                
                \item[Case 2.] $X_\alpha^\geq =X_\alpha\cup X_{\alpha+2}$ with $X_{\alpha+2}=\{a\}$.
                
                Since for every $b\in X_\alpha$ we have $\{a,b\}\in D$, then $|X_\alpha^\geq|\leq \alpha+3$, and we have
                \begin{align*}
                    2|D| &= \sum_{x\in [n]} i_x =i_a+ \sum_{x\in X_\alpha} i_x + \sum_{x\in X_{\alpha-1}^\leq} i_x\leq\\
                    &\leq (\alpha +2)+\underbrace{|X_{\alpha}|}_{|X_{\alpha}^\geq|-1}\alpha+\underbrace{|X_{\alpha-1}^\leq|}_{n-|X_\alpha^\geq|}(\alpha -1)=\\
                    &=n(\alpha-1)+|X_\alpha^\geq|+2\leq n\alpha -n + \alpha+5.
                \end{align*}
                Since $2|D|>\alpha n$, it turns out that $n<\alpha+5$. By hypothesis, $n\geq 2\alpha +3$, then $2\alpha+3<\alpha +5$, that implies $\alpha < 2$, which cannot be true since $\alpha\geq 2$.
            \end{description}
            In both cases, we arise to a contradiction. Therefore, we can conclude that $\alpha n\geq 2k+4\alpha$ as claimed.
        \end{proofclaim}

        \begin{claim}
            If $\Kneser{n}{2}$ admits a $k$-tuple dominating set with cardinality $k+2\alpha-1$, then $(\alpha-1) n\geq 2k+3(\alpha-1)$.
        \end{claim}
        
            \begin{proofclaim}
            Let $D$ be a $k$-tuple dominating set of $\Kneser{n}{2}$ with cardinality $|D|=k+2\alpha-1$. We are intended to prove that $(\alpha-1)n\geq 2k+3(\alpha-1)$ or equivalently $(\alpha-1) n+\alpha+1\geq 2|D|$. Let us suppose, on the contrary, that $(\alpha-1) n+ \alpha +1 < 2|D|$.
            
            By Lemma \ref{lemma: cotas indices} for every pair of elements $a,b\in[n]$ we have
            \begin{equation}
            \label{eq: kupla general rayo chico}
                i_a+i_b\leq\begin{cases}
                |D|-k+2=2\alpha+ 1 & \text{if }\{a,b\}\in D,\\
                |D|-k=2\alpha -1 & \text{if }\{a,b\}\notin D.
                \end{cases}
            \end{equation}
            Since $(\alpha-1) n < 2|D|=\sum_{x\in [n]}i_x$, we have that there is at least one element $a\in [n]$ for which $i_a\geq \alpha$. For every other element $b\in X_{\alpha}^\geq$, it holds $i_a+i_b\geq 2\alpha$, and by (\ref{eq: kupla general rayo chico}) $\{a,b\}\in D$.
            
            Let us see that $|X_{\alpha}^\geq|>1$. On the contrary, $X_{\alpha}^\geq=\{a\}$ and
            \begin{align*}
                (\alpha-1) n +\alpha +1 < 2|D| &= \sum_{x\in [n]} i_x = i_a + \sum_{x\neq a} \underbrace{i_x}_{\leq \alpha-1} \leq i_a + (\alpha-1) n  - (\alpha-1).
            \end{align*}
            Thus, $i_a > 2\alpha $, which cannot be true since $i_x\leq 2\alpha$ for every element $x\in [n]$ by Remark \ref{rem: indices}. Therefore, $|X_{\alpha}^\geq|\geq 2$.
            
            Then, we have that for every $x\in X_{\alpha}^\geq$, $i_x\leq \alpha +1$. In fact, if for some element $y$ it holds that $i_y>\alpha +1$, then by (\ref{eq: kupla general rayo chico}) for every other $x\in [n]$
            $$i_x\leq 2\alpha +1-i_y<2\alpha+1-(\alpha+1)=\alpha,$$
            and $y$ would be the only element in $X_{\alpha}^\geq$, but $|X_{\alpha}^\geq|>1$. 
            Therefore, $i_x\leq \alpha +1$ for every $x\in X_{\alpha}^\geq$. What is more, there is at most one element in $X_{\alpha +1}$ since if $x,y\in X_{\alpha +1}$, then $i_x+i_y=2\alpha+2$ and this cannot be true by (\ref{eq: kupla general rayo chico}). 
            
            Therefore, either $X_{\alpha}^\geq=X_{\alpha}$ or $X_{\alpha}^\geq=X_{\alpha}\cup X_{\alpha+1}$ with $|X_{\alpha+1}|=1$.
            \begin{description}
                \item[Case 1.] $X_{\alpha}^\geq =X_{\alpha}$.
                
                For every pair $a,b\in X_{\alpha}$ we have $\{a,b\}\in D$ since $i_a+i_b=2\alpha$ and it holds (\ref{eq: kupla general rayo chico}). Thus $\binom{X_{\alpha}}{2}\subseteq D$ and we have $|X_{\alpha}|\leq \alpha+1$. Then
                \begin{align*}
                    2|D| &= \sum_{x\in [n]} i_x = \sum_{x\in X_{\alpha}} i_x + \sum_{x\in X_{\alpha-1}^\leq} i_x\leq\\
                    &\leq |X_{\alpha}|\alpha+\underbrace{|X_{\alpha-1}^\leq|}_{n-|X_{\alpha}|}(\alpha-1)=\\
                    &=n(\alpha-1)+|X_{\alpha}|\leq n(\alpha-1) +\alpha+1.
                \end{align*}
                But we have assumed that $2|D|>n(\alpha-1) +\alpha+1$.
                
                \item[Case 2.] $X_{\alpha}^\geq =X_{\alpha}\cup X_{\alpha+1}$ with $X_{\alpha+1}=\{a\}$.
                
                For every other $b\in X_{\alpha-1}^\geq$ we have $i_a+i_b\geq 2\alpha$ thus $\{a,b\}\in D$ and in consequence $|X_{\alpha-1}^\geq|\leq \alpha+2$. If $X=X_{\alpha-1}\cup X_{\alpha}$, we have
                \begin{align*}
                    2|D| &= \sum_{x\in [n]} i_x =i_a+ \sum_{x\in X} i_x + \sum_{x\in X_{\alpha-2}^\leq} i_x\leq\\
                    &\leq (\alpha +1)+\underbrace{|X|}_{|X_{\alpha-1}^\geq|-1}\alpha+\underbrace{|X_{\alpha-2}^\leq|}_{n-|X_{\alpha-1}^\geq|}(\alpha -2)=\\
                    &=n(\alpha-2)+2|X_{\alpha-1}^\geq|+1\leq n(\alpha-1) -n + 2\alpha+5.
                \end{align*}
                Since $2|D|>(\alpha-1) n+\alpha +1$, it turns out that $n<\alpha+4$, but by hypothesis $n\geq 2\alpha+3\geq \alpha+5$, since $\alpha\geq 2$.
            \end{description}
            In both cases, we arise to a contradiction. Therefore, we can conclude that $(\alpha-1) n\geq 2k+3(\alpha-1)$ as claimed.
        \end{proofclaim}

        Assume that $n=\left\lceil\frac{2}{\alpha}k\right\rceil+3$. Suppose that $K(n,2)$ admits a $k$-tuple dominating set with cardinality $k+2\alpha$. From the first claim, it follows that $\alpha n\geq 2k+4\alpha$, or, equivalently,
        $$n\geq \frac{2}{\alpha}k+4,$$
        which does not hold. Thus, there does not exist any $k$-tuple dominating set with $k+2\alpha$ vertices and
        $$\kdom{n,2}{k}\geq k+2\alpha+1.$$
        
        Similarly, if $\frac{2}{\alpha}k+4\leq n<\frac{2}{\alpha-1}k+3$ and $K(n,2)$ admits a $k$-tuple dominating set with cardinality $k+2\alpha-1$, from the second claim it follows that $(\alpha-1)n\geq 2k+3(\alpha-1)$, i.e.
        $$n\geq \frac{2}{\alpha-1}k+3,$$
        which is not true. Then, there does not exist any $k$-tuple dominating set with $k+2\alpha-1$ vertices and
        $$\kdom{n,2}{k}\geq k+2\alpha.$$
    \end{proof}
\end{proposition}

\subsection*{\texorpdfstring{Construction of $k$-tuple dominating sets for $\Kneser{n}{2}$}{Construction of k-tuple dominating sets for K(n,2)}}

In order to provide $k$-tuple dominating sets of $\Kneser{n}{2}$ whose cardinalities achieve the lower bounds for $\kdom{n,2}{k}$ given in Proposition \ref{prop: lower bounds kdom(n,2) V2}, let us introduce the following definition.

\begin{definition}\label{def: conjuntos D(h,alfa)}
    Let $m\in\mathbb{N}$. For $i\in \mathbb{N}$, $i<\frac{m}{2}$, we define $D_i^{[m]}$ as the set 
    $$D_i^{[m]}=\big\{\{\xi,\xi+i\}:\xi\in[m]\big\},$$
    where the sums are taken modulo $m$. And for $\alpha\in\mathbb{N}$, $\alpha\geq 2$, such that $m>\alpha$, let $a=\left\lfloor\frac{\alpha}{2}\right\rfloor$. We define $D^{m,\alpha}$ as the set given by
    $$D^{m,\alpha} = \left( \bigcup_{i=1}^{a} D_i^{[m]} \right)\cup \widehat{D},$$
    with
    $$\widehat{D}=\begin{cases}
    \left\{\left\{\xi,\xi+\left\lfloor\frac{m}{2}\right\rfloor\right\},\xi\in \left[\left\lfloor\frac{m}{2}\right\rfloor\right]\right\}, &\text{if $\alpha$ is odd,}\\
    \emptyset,&\text{if $\alpha$ is even.}
\end{cases}$$
\end{definition}

\begin{example}\label{ej: ejemplos D(m,alfa)}
    \begin{itemize}

        \item Let $\alpha=3$, $m=6>\alpha$ and  $a=\left\lfloor\frac{\alpha}{2}\right\rfloor=1$. The set $D^{m,\alpha}=D^{6,3}$ is given by
        $$D^{6,3}=D_1^{[6]}\cup \widehat{D},$$
        where
        \begin{align*}
            D_1^{[6]} &= \big\{\{1,2\},\{2,3\},\{3,4\},\{4,5\},\{5,6\},\{6,1\}\big\},\\
            \widehat{D} &= \big\{\{1,4\},\{2,5\},\{3,6\}\big\}.
        \end{align*}
        $D^{6,3}$ is a set of vertices of $\Kneser{n}{2}$ for every $n\geq 6$ and the occurrences of the elements in $[n]$ for the set $D^{6,3}$ are
        $$i_x=\begin{cases}
            3=\alpha, &\text{if }x\in [6],\\
            0, &\text{otherwise}.
        \end{cases}$$
        
        \begin{figure}[ht]
            \centering
            \begin{tikzpicture}[rotate=-15]
            
                \foreach \x in {0,30,...,330} {
                \draw (\x:3) -- (\x+30:3);
                }
                \foreach \x in {0,30,...,150} {
                \draw (\x:3) -- (\x+180:3);
                \draw ({2*\x}:3) to[bend left=20] ({2*\x+60}:3);
                \draw ({2*\x+30}:3) to[bend right=20] ({2*\x+150}:3);
                }

                \foreach \x in {45,165,285}{
                \draw ({\x-45}:3) -- (\x:1.5);
                \draw ({\x+45}:3) -- (\x:1.5);
                \draw ({\x+135}:3) to[bend left = 20] (\x:1.5);
                \draw ({\x-135}:3) to[bend right = 20] (\x:1.5);
                }

                \draw (45:1.5) -- (165:1.5) -- (285:1.5) -- cycle;

                \filldraw (45:1.5) circle (2.5pt);
                \filldraw (165:1.5) circle (2.5pt);
                \filldraw (285:1.5) circle (2.5pt);
                
                \foreach \x in {0,30,...,330} {
                \filldraw (\x:3) circle (2.5pt);
                }
                
                \foreach \x in {30,90,...,330} {
                \filldraw[fill=white] (\x:3) circle (2.5pt);
                }
                \foreach \x in {45,165,285} {
                \filldraw[fill=white] (\x:1.5) circle (2.5pt);
                }

                \node (35) at (0:3) {};
                \node (12) at (30:3) {};
                \node (46) at (60:3) {};
                \node (23) at (90:3) {};
                \node (15) at (120:3) {};
                \node (34) at (150:3) {};
                \node (26) at (180:3) {};
                \node (45) at (210:3) {};
                \node (13) at (240:3) {};
                \node (56) at (270:3) {};
                \node (24) at (300:3) {};
                \node (16) at (330:3) {};
                \node (14) at (45:1.5) {};
                \node (36) at (165:1.5) {};
                \node (25) at (285:1.5) {};

                \foreach \x in {12,13,14,15,16,23,24,25,26,34,35,36,45,46,56} {
                \begin{scope}[shift={(\x)}, rotate=105]
                    \filldraw [fill=white, thick] (0,0) circle (0.2);
                    \draw (0:0.2) -- (180:0.2) (60:0.2) -- (240:0.2) (120:0.2) -- (300:0.2);
                \end{scope}}

            \begin{scope}[shift={(0:3)}, rotate=105]
                \filldraw[fill = gray!50]  (0,0) -- (-120:0.2) arc(-120:-180:0.2) -- cycle;
                \filldraw[fill = gray!50]  (0,0) -- (-240:0.2) arc(-240:-300:0.2) -- cycle;
            \end{scope}

            \begin{scope}[shift={(30:3)}, rotate=105]
                \filldraw[fill = red]  (0,0) -- (0:0.2) arc(0:-60:0.2) -- cycle;
                \filldraw[fill = red]  (0,0) -- (-60:0.2) arc(-60:-120:0.2) -- cycle;
            \end{scope}

            \begin{scope}[shift={(60:3)}, rotate=105]
                \filldraw[fill = gray!50]  (0,0) -- (-180:0.2) arc(-180:-240:0.2) -- cycle;
                \filldraw[fill = gray!50]  (0,0) -- (-300:0.2) arc(-300:-360:0.2) -- cycle;
            \end{scope}

            \begin{scope}[shift={(90:3)}, rotate=105]
                \filldraw[fill = red]  (0,0) -- (-60:0.2) arc(-60:-120:0.2) -- cycle;
                \filldraw[fill = red]  (0,0) -- (-120:0.2) arc(-120:-180:0.2) -- cycle;
            \end{scope}

            \begin{scope}[shift={(120:3)}, rotate=105]
                \filldraw[fill = gray!50]  (0,0) -- (0:0.2) arc(0:-60:0.2) -- cycle;
                \filldraw[fill = gray!50]  (0,0) -- (-240:0.2) arc(-240:-300:0.2) -- cycle;
            \end{scope}

            \begin{scope}[shift={(150:3)}, rotate=105]
                \filldraw[fill = red]  (0,0) -- (-120:0.2) arc(-120:-180:0.2) -- cycle;
                \filldraw[fill = red]  (0,0) -- (-180:0.2) arc(-180:-240:0.2) -- cycle;
            \end{scope}

            \begin{scope}[shift={(180:3)}, rotate=105]
                \filldraw[fill = gray!50]  (0,0) -- (-60:0.2) arc(-60:-120:0.2) -- cycle;
                \filldraw[fill = gray!50]  (0,0) -- (-300:0.2) arc(-300:-360:0.2) -- cycle;
            \end{scope}

            \begin{scope}[shift={(210:3)}, rotate=105]
                \filldraw[fill = red]  (0,0) -- (-180:0.2) arc(-180:-240:0.2) -- cycle;
                \filldraw[fill = red]  (0,0) -- (-240:0.2) arc(-240:-300:0.2) -- cycle;
            \end{scope}

            \begin{scope}[shift={(240:3)}, rotate=105]
                \filldraw[fill = gray!50]  (0,0) -- (0:0.2) arc(0:-60:0.2) -- cycle;
                \filldraw[fill = gray!50]  (0,0) -- (-120:0.2) arc(-120:-180:0.2) -- cycle;
            \end{scope}

            \begin{scope}[shift={(270:3)}, rotate=105]
                \filldraw[fill = red]  (0,0) -- (-240:0.2) arc(-240:-300:0.2) -- cycle;
                \filldraw[fill = red]  (0,0) -- (-300:0.2) arc(-300:-360:0.2) -- cycle;
            \end{scope}

            \begin{scope}[shift={(300:3)}, rotate=105]
                \filldraw[fill = gray!50]  (0,0) -- (-60:0.2) arc(-60:-120:0.2) -- cycle;
                \filldraw[fill = gray!50]  (0,0) -- (-180:0.2) arc(-180:-240:0.2) -- cycle;
            \end{scope}

            \begin{scope}[shift={(330:3)}, rotate=105]
                \filldraw[fill = red]  (0,0) -- (0:0.2) arc(0:-60:0.2) -- cycle;
                \filldraw[fill = red]  (0,0) -- (-300:0.2) arc(-300:-360:0.2) -- cycle;
            \end{scope}

            \begin{scope}[shift={(45:1.5)}, rotate=105]
                \filldraw[fill = red]  (0,0) -- (0:0.2) arc(0:-60:0.2) -- cycle;
                \filldraw[fill = red]  (0,0) -- (-180:0.2) arc(-180:-240:0.2) -- cycle;
            \end{scope}

            \begin{scope}[shift={(165:1.5)}, rotate=105]
                \filldraw[fill = red]  (0,0) -- (-120:0.2) arc(-120:-180:0.2) -- cycle;
                \filldraw[fill = red]  (0,0) -- (-300:0.2) arc(-300:-360:0.2) -- cycle;
            \end{scope}

            \begin{scope}[shift={(285:1.5)}, rotate=105]
                \filldraw[fill = red]  (0,0) -- (-60:0.2) arc(-60:-120:0.2) -- cycle;
                \filldraw[fill = red]  (0,0) -- (-240:0.2) arc(-240:-300:0.2) -- cycle;
            \end{scope}

            \end{tikzpicture}
            \caption{Set $D^{6,3}$ in the Kneser graph $\Kneser{6}{2}$}
            \label{fig: D(6,3) en Kneser(6,2)}
        \end{figure}
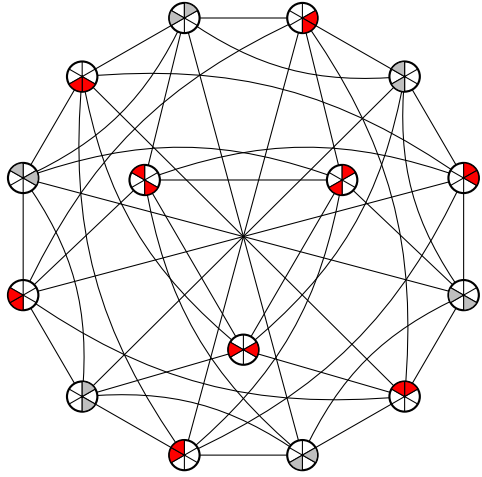

        \item Let $\alpha=4$, $m=14>\alpha$ and $a=\left\lfloor\frac{\alpha}{2}\right\rfloor=2$. The set $D^{m,\alpha}=D^{14,4}$ is given by
        $$D^{14,4}=D_1^{[14]}\cup D_{2}^{[14]},$$
        where
        \begin{align*}
            D_1^{[14]} = \big\{&\{1,2\},\{2,3\},\{3,4\},\{4,5\},\{5,6\},\{6,7\},\{7,8\},\{8,9\},\\
            &\{9,10\},\{10,11\},\{11,12\},\{12,13\},\{13,14\},\{14,1\}\big\},\\
            D_2^{[14]} = \big\{&\{1,3\},\{2,4\},\{3,5\},\{4,6\},\{5,7\},\{6,8\},\{7,9\},\{8,10\},\\
            &\{9,11\},\{10,12\},\{11,13\},\{12,14\},\{13,1\},\{14,2\}\big\}.
        \end{align*}
        $D^{14,4}$ is a set of vertices of $\Kneser{n}{2}$ for every $n\geq 14$ and the occurrences of the elements in $[n]$ for the set $D^{14,4}$ are
        $$i_x=\begin{cases}
            4=\alpha, &\text{if }x\in [14],\\
            0, &\text{otherwise}.
        \end{cases}$$

        \item Let $\alpha=3$, $m=11>\alpha$, and $a=\left\lfloor\frac{\alpha}{2}\right\rfloor=1$. The set $D^{m,\alpha}=D^{11,3}$ is given by
        $$D^{11,3}=D_1^{[11]}\cup \widehat{D},$$
        where
        \begin{align*}
            D_1^{[11]} &= \big\{\{1,2\},\{2,3\},\{3,4\},\{4,5\},\{5,6\},\{6,7\},\{7,8\},\{8,9\},\{9,10\},\{10,11\},\{11,1\}\big\},\\
            \widehat{D} &= \big\{\{1,6\},\{2,7\},\{3,8\},\{4,9\},\{5,10\}\big\}.
        \end{align*}
        $D^{11,3}$ is a set of vertices of $\Kneser{n}{2}$ for every $n\geq 11$ and the occurrences of the elements in $[n]$ for the set $D^{11,3}$ are
        $$i_x=\begin{cases}
            3=\alpha, &\text{if }x\in [10],\\
            2=\alpha-1, &\text{if }x=11,\\
            0, &\text{otherwise}.
        \end{cases}$$
    \end{itemize}
\end{example}

Let us observe that $D^{m,\alpha}\subseteq\binom{[m]}{2}$. In consequence, $D^{m,\alpha}$ is a set of vertices of $\Kneser{n}{2}$ for every $n\geq m$, and any element in $[m]$ has at most $\alpha$ occurrences in $D^{m,\alpha}$. In fact, $m\geq\alpha+1>2a$, i.e., $a<\frac{m}{2}$. So, each element $x\in [m]$ appears in exactly two vertices of $D_i^{[m]}$ for each $1\leq i\leq a$. Moreover, for $i\neq j$, $D_i^{[m]}\cap D_j^{[m]}=\emptyset$ since $i+j\leq2a<m$. On the other hand, for $\alpha$ odd, a vertex in $\widehat{D}$ is not in $D_i^{[m]}$ for any $i$. Furthermore, every element $x\in [m]$ appears in exactly one vertex of $\widehat{D}$ if $m$ is even, and every element $x\in [m-1]$ appears in exactly one vertex of $\widehat{D}$ if $m$ is odd.

Thus, for both $\alpha$ and $m$ odd we have
\begin{equation}\label{eq: ocurrencias D(h,alfa) caso impar}
    i_x\left(D^{m,\alpha}\right)=\begin{cases}
    \alpha, &\text{if }x\in [m-1],\\
    \alpha-1=2a, &\text{if }x=m,\\
    0, &\text{otherwise}.
\end{cases}
\end{equation}
and if either $\alpha$ or $m$ are even, then
\begin{equation}\label{eq: ocurrencias D(h,alfa) caso par}
    i_x\left(D^{m,\alpha}\right)=\begin{cases}
    \alpha, &\text{if }x\in [m],\\
    0, &\text{otherwise}.
\end{cases}
\end{equation}
In any case, we have 
\begin{equation}\label{eq: cardinal D(h,alfa)}
    \left|D^{m,\alpha}\right|=\frac{1}{2}\sum_{x\in [n]}i_x=\left\lfloor\frac{\alpha m}{2}\right\rfloor.
\end{equation}

\begin{lemma}\label{lem: kupla set cardinality k+2 alfa}
    Let $\alpha,k\in\mathbb{N}$ with $\alpha\geq 2$ and $n\geq \alpha+2$. If $n \geq \frac{2}{\alpha}k+4$, then there exists a $k$-tuple dominating set of $\Kneser{n}{2}$ with cardinality $k+2\alpha$.

    \begin{proof}
    Due to Lemma \ref{lemma: cotas indices} a set $D$ is a $k$-tuple dominating set if and only if it verifies, for every pair of elements $x$ and $y$ in $[n]$
        $$i_x+i_y\leq 
        \begin{cases}
            |D|-k+2=2\alpha+2, & \text{if }\{x,y\}\in D,\\
            |D|-k=2\alpha, & \text{if }\{x,y\}\notin D.
        \end{cases}$$
        Therefore, it is enough to find a set $D$ of vertices with cardinality $k+2\alpha$ such that $i_x\leq \alpha$ for each $x\in [n]$. Let $a=\left\lfloor\frac{\alpha}{2}\right\rfloor$. Since $n\geq\alpha+2>2a+1$, we have $a<\left\lfloor\frac{n}{2}\right\rfloor$. Thus, let us consider the set $D^{n,\alpha}$ as in Definition \ref{def: conjuntos D(h,alfa)}. Note that $\left|D^{n,\alpha}\right|=\left\lfloor\frac{\alpha n}{2}\right\rfloor$ by (\ref{eq: cardinal D(h,alfa)}).
        Moreover, $|D^{n,\alpha}|\geq k+2\alpha$. In fact, as $\frac{2}{\alpha}k+4\leq n$, we have $\alpha n\geq 2(k+2\alpha)$. If either $\alpha$ or $n$ are even, then $$2|D^{n,\alpha}|=n\alpha \geq 2(k+2\alpha).$$
        On the other hand, if both $\alpha$ and $n$ are odd then
        $$2|D^{n,\alpha}|= n\alpha-1\geq 2(k+2\alpha)-1.$$
        As $2|D^{n,\alpha}|$ is an even integer and $2(k+2\alpha)-1$ is an odd integer, it turns out that $2|D^{n,\alpha}|>2(k+2\alpha)-1$ and in consequence $2|D^{n,\alpha}|\geq 2(k+2\alpha)$. 

        In any case, $|D^{n,\alpha}|\geq k+2\alpha$ as claimed. 
        
        Note that $i_x(D^{n,\alpha})\leq \alpha$ for every $x\in [n]$, by (\ref{eq: ocurrencias D(h,alfa) caso impar}) and (\ref{eq: ocurrencias D(h,alfa) caso par}). Thus, eliminating any $|D^{n,\alpha}|-(k+2\alpha)$ vertices from $D^{n,\alpha}$ gives as a result a set of vertices $D$ with cardinality $k+2\alpha$ such that $i_x(D)\leq i_x(D^{n,\alpha})\leq \alpha$ for every $x\in [n]$. This is, a $k$-tuple dominating set with cardinality $k+2\alpha$.
        
        Therefore, $\Kneser{n}{2}$ admits a $k$-tuple dominating set with cardinality $k+2\alpha$.
    \end{proof}
\end{lemma}

\begin{example}
\begin{itemize}
    \item Let $\alpha=3$, $n=11$ and $k=9$. We have
    $\frac{2}{\alpha}k+4\leq n < \frac{2}{\alpha -1}k+3$. We are in the hypothesis of Proposition \ref{lem: kupla set cardinality k+2 alfa}. Let us consider the set $D^{11,3}$ given in Example \ref{ej: ejemplos D(m,alfa)}. We have $\left|D^{11,3}\right|=16>15=k+2\alpha$. Thus, by eliminating any vertex of $D^{11,3}$ we get a $k$-tuple dominating set $D$ of $\Kneser{11}{2}$ with cardinality $k+2\alpha$.
    
    \item Let $\alpha=4$, $n=14$ and $k=17$. We have
    $\frac{2}{\alpha}k+4\leq n < \frac{2}{\alpha -1}k+3$. We are in the hypothesis of Proposition \ref{lem: kupla set cardinality k+2 alfa}. Let us consider the set $D^{14,4}$ given in Example \ref{ej: ejemplos D(m,alfa)}. We have $\left|D^{14,4}\right|=28>25=k+2\alpha$. Thus, by eliminating any three vertices of $D^{14,4}$ we get a $k$-tuple dominating set $D$ of $\Kneser{14}{2}$ with cardinality $k+2\alpha$.
\end{itemize}

\end{example}

\begin{lemma}\label{lem: kupla set cardinality k+2 alfa +1}
    Let $\alpha,k\in\mathbb{N}$ with $\alpha\geq 2$ and $n\geq 2\alpha+3$. If $n =\left\lceil\frac{2}{\alpha}k\right\rceil+3$, then there exists a $k$-tuple dominating set of $\Kneser{n}{2}$ with cardinality $k+2\alpha+1$.

    \begin{proof}
        Due to Lemma \ref{lemma: cotas indices} a set $D$ is a $k$-tuple dominating set if and only if for every pair of elements $x$ and $y$ in $[n]$ it holds
        $$i_x+i_y\leq 
        \begin{cases}
            |D|-k+2=2\alpha+3 & \text{if }\{x,y\}\in D,\\
            |D|-k=2\alpha+1 & \text{if }\{x,y\}\notin D.
        \end{cases}$$
        Therefore it is enough to find a set $D$ of vertices with cardinality $k+2\alpha+1$ such that $\alpha\leq i_x\leq \alpha+1$ for each $x\in [n]$, and $\{x,y\}\in D$ for every pair $x,y\in X_{\alpha+1}$.

        Let $\lambda\in \mathbb{N}$ such that $k=\alpha\lambda+b$ with $0\leq b\leq \alpha-1$. And let $a=\left\lfloor \frac{\alpha}{2}\right\rfloor$. We will give such a set $D$ in terms of $b$ considering the cases $b=0$, $1\leq b \leq a$ and $a+1\leq b \leq \alpha-1$. 

        In order to do so, let us consider for $h\in\mathbb{N}$, $1\leq h< n-\alpha$ the set
        \begin{equation}\label{eq: definicion D(h)}
            D(h)=D^{n-h,\alpha}\cup \binom{[n-h+1..n]}{2}.
        \end{equation}
        Let us note that if $(n-h)\alpha$ is even, then using (\ref{eq: ocurrencias D(h,alfa) caso par}) we have that each element $x\in [n-h]$ has exactly $\alpha$ occurrences in $D^{n-h,\alpha}$ and none in $\binom{[n-h+1..n]}{2}$, and each $x\in [n-h+1..n]$ has exactly $h-1$ occurrences in $\binom{[n-h+1..n]}{2}$ and none in $D^{n-h,\alpha}$. In consequence
        \begin{equation}\label{eq: ocurrencias D(h)}
            i_x(D(h))=\begin{cases}
                \alpha, &\text{if }x\in [n-h],\\
                h-1, &\text{if }x\in [n-h+1..n].
            \end{cases}
        \end{equation}
        Furthermore, by (\ref{eq: cardinal D(h,alfa)}) we have
        \begin{equation}\label{eq: cardinal D(h)}
            \left|D(h)\right|=\left|D^{n-h,\alpha}\right|+\left|\binom{[n-h+1..n]}{2}\right|=\frac{(n-h)\alpha}{2}+\binom{h}{2}.
        \end{equation}
        
        In the first case, $b=0$, we will prove that for $h=\alpha+2$ the set $D(h)$ is itself a $k$-tuple dominating set with cardinality $k+2\alpha$. In the remaining cases, we will consider the set $D(h)$ for $h=2b+2$ and $h=2b-\alpha+2$ respectively and modify them with the aim of giving $k$-tuple dominating sets of the desired cardinality.

        \begin{description}
            \item[Case 1.] $b=0$.

            Let $h=\alpha+2$. Since $n\geq 2\alpha+3$, we have $n-h\geq \alpha+1>\alpha$. So, we consider the set $D=D(h)$ as in (\ref{eq: definicion D(h)}). Note that $n=\left\lceil\frac{2}{\alpha}k\right\rceil+3=2\lambda+3$. Thus, for $\alpha$ odd, we have both $n$ and $h$ are odd and in consequence $n-h$ is even. Thus, $(n-h)\alpha$ is even. From (\ref{eq: ocurrencias D(h)}) it turns out that $\alpha\leq i_x\leq \alpha+1$ for each $x\in [n]$, and $\{x,y\}\in D$ for every pair $x,y\in X_{\alpha+1}$. Thus, $D$ is a $k$-tuple dominating set. What is more,
            $$2|D|=\sum_{x\in [n]}i_x=(n-h)\alpha+h(\alpha+1)=\alpha n+h=\alpha \left(\frac{2}{\alpha}k+3\right)+\alpha+2=2(k+2\alpha+1).$$
            Therefore, $|D|$ is a $k$-tuple dominating set with cardinality $k+2\alpha+1$.

            \item[Case 2.] $1\leq b\leq a$.

            Let $h=2b+2$. $h$ is an even integer between $4$ and $2a+2\leq\alpha+2$. Since $n\geq 2\alpha+3$, we have $n-h\geq (2\alpha+3)-(\alpha+2)=\alpha+1>\alpha$. So, we consider the set $D(h)$ as in (\ref{eq: definicion D(h)}).

            Note that $n=\left\lceil\frac{2}{\alpha}k\right\rceil+3=2\lambda+4$. Thus, both $n$ and $h$ are even and in consequence also $n-h$. From (\ref{eq: cardinal D(h)}),
            \begin{equation}
                \label{eq: cardinal D(h) Caso 2}
                2|D(h)|=(n-h)\alpha+h(h-1)=(n\alpha+h)-h\underbrace{(\alpha-h+2)}_{\alpha-2b}=2(k+2\alpha+1)-\underbrace{h(\alpha-2b)}_{\geq 0},
            \end{equation}
            since 
            $$2(k+2\alpha+1)=2k+4\alpha+2=2(\lambda\alpha+b)+4\alpha+2=(2\lambda+4)\alpha+(2b+2)=n\alpha+h.$$

            If $b<a$, let us consider the following set $D$.
            \begin{equation}\notag
                D = \Bigg[ D(h)\setminus\underbrace{\Bigg(\smashoperator[r]{\bigcup_{\substack{1\leq j \leq b+1\\1\leq\xi\leq\alpha-2b}}}\ \big\{\{\xi,\xi+j\}\big\}\Bigg)}_{D_1}\Bigg]\cup\underbrace{\Bigg(\smashoperator[r]{\bigcup_{\substack{1\leq j \leq b+1\\1\leq\xi\leq\alpha-2b}}}\ \big\{\{\xi,n-h+2j-1\},\{\xi+j,n-h+2j\}\big\}\Bigg)}_{D_2},
            \end{equation}
            where the sums in $D_1$ are taken modulo $n-h$.
            
            Let us see that $D_1\subseteq D(h)$. In fact, we have that $b+1\leq a$, so for each $j$ it holds $\{\xi,\xi+j\}\in D_j^{[n-h]}\subseteq D(h)$. What is more, $\alpha-2b<\alpha<n-h$, so for fixed $j$, the vertices $\{\xi,\xi+j\}$ are different. In consequence, $|D_1|=(b+1)(\alpha-2b)$.
            
            On the other hand, $D_2\cap D(h)=\emptyset$ since no vertex in $D(h)$ has an element in $[n-h]$ and other in $[n-h+1..n]$. And the vertices in $D_2$ are different. Thus, $|D_2|=2(b+1)(\alpha-2b)$.

            Note that for each $x\in [n]$, we have
            $$i_x(D)=i_x(D(h))-i_x(D_1)+i_x(D_2).$$
            For $x\in [n-h]$, we have $i_x(D_1)=i_x(D_2)$, and in consequence 
            $$i_x(D)=i_x(D(h))=\alpha,$$
            and for $x\in [n-h+1..n]$, we have $i_x(D_1)=0$ and $i_x(D_2)=\alpha-2b$, and in consequence 
            $$i_x(D)=\underbrace{i_x(D(h))}_{h-1}+\alpha-2b=\alpha+1.$$
            We have $\alpha\leq i_x(D)\leq \alpha+1$ for each $x\in [n]$, and $\{x,y\}\in D$ for every pair $x,y\in X_{\alpha+1}$. Thus, $D$ is a $k$-tuple dominating set. Furthermore,
            $$|D|=|D(h)|-|D_1|+|D_2|=k+2\alpha+1.$$
            Therefore, $|D|$ is a $k$-tuple dominating set with cardinality $k+2\alpha+1$.

            If $b=a$ and $\alpha$ is even then the set $D(h)$ is itself a $k$-tuple dominating set with cardinality $k+2\alpha+1$ since (\ref{eq: cardinal D(h) Caso 2}) holds. If $\alpha$ is odd, then let us consider the following set $D$.
            \begin{equation}\notag
                D = \Bigg[ D(h)\setminus\underbrace{\Bigg(\smashoperator[r]{\bigcup_{{1\leq\xi\leq b+1}}}\ \big\{\{\xi,\xi+1\}\big\}\Bigg)}_{D_1}\Bigg]\cup\underbrace{\Bigg(\smashoperator[r]{\bigcup_{1\leq\xi\leq b+1}}\ \big\{\{\xi,n-h+2\xi-1\},\{\xi+1,n-h+2\xi\}\big\}\Bigg)}_{D_2}.
            \end{equation}
            Note that $D_1\subseteq D_1^{[n-h]}\subseteq D(h)$. Moreover, since $\alpha$ is odd, we have $\alpha=2a+1=2b+1$. Thus, $n-h\geq (2\alpha+3)-(2b+2)=\alpha+2>b+1$. So, the vertices in $D_1$ are different and $|D_1|=b+1$. As in the case $b<a$, $D_2\cap D(h)=\emptyset$ and the vertices in $D_2$ are different. Thus, $|D_2|=2b+2$. We have that for each $x\in [n]$ it holds
            $$i_x(D)=i_x(D(h))-i_x(D_1)+i_x(D_2).$$
            We have $i_x(D_1)=i_x(D_2)$ if $x\in [n-h]$ and $i_x(D_1)=0$, $i_x(D_2)=1$ if $x\in[n-h+1..n]$. In consequence
            $$i_x(D)=\begin{cases}
                \alpha, &\text{if }x\in [n-h],\\
                h=\alpha+1, &\text{if }x\in [n-h+1..n].
            \end{cases}$$
            This is, $\alpha\leq i_x(D)\leq \alpha+1$ for each $x\in [n]$, and $\{x,y\}\in D$ for every pair $x,y\in X_{\alpha+1}$. $D$ is a $k$-tuple dominating set with cardinality
            $$|D|=|D(h)|-|D_1|+|D_2|=k+2\alpha+1.$$
            
            \item[Case 3.] $a+1 \leq b\leq\alpha-1$.

            Let $h=2b-\alpha+2$. We have $3\leq h\leq \alpha$, and $h$ has the same parity as $\alpha$, and $n-h\geq (2\alpha+3)-\alpha=\alpha+3>\alpha$. So, we consider the set $D(h)$ as in (\ref{eq: definicion D(h)}). 
            
            Note that $n=\left\lceil\frac{2}{\alpha}k\right\rceil+3=2\lambda+5$. Besides, if $\alpha$ is odd, then both $n$ and $h$ are odd and in consequence $n-h$ is even. So, in any case, $(n-h)\alpha$ is even. From (\ref{eq: cardinal D(h)}),
            \begin{equation}
                \label{eq: cardinal D(h) Caso 3}
                2|D(h)|=(n-h)\alpha+h(h-1)=(n\alpha+h)-h\underbrace{(\alpha-h+2)}_{2\alpha-2b}=2(k+2\alpha+1)-\underbrace{2h(\alpha-b)}_{> 0},
            \end{equation}
            since 
            $$2(k+2\alpha+1)=2k+4\alpha+2=2(\lambda\alpha+b)+4\alpha+2=(2\lambda+5)\alpha+(2b-\alpha+2)=n\alpha+h.$$
            
            The procedure is similar than Case 2. This is, we eliminate $h(\alpha-b)$ vertices from $D(h)$ and add another $2h(\alpha-b)$ vertices in order to get a $k$-tuple dominating set with cardinality $k+2\alpha+1$.
            
            If $\alpha$ is even, let us consider the following set $D$.
            \begin{equation}\notag
                D = \Bigg[ D(h)\setminus\underbrace{\Bigg(\smashoperator[r]{\bigcup_{\substack{1\leq j \leq b-a+1\\1\leq\xi\leq2\alpha-2b}}}\ \big\{\{\xi,\xi+j\}\big\}\Bigg)}_{D_1}\Bigg]\cup\underbrace{\Bigg(\smashoperator[r]{\bigcup_{\substack{1\leq j \leq b-a+1\\1\leq\xi\leq2\alpha-2b}}}\ \big\{\{\xi,n-h+2j-1\},\{\xi+j,n-h+2j\}\big\}\Bigg)}_{D_2},
            \end{equation}
            where the sums in $D_1$ are taken modulo $n-h$.

            Let us see that $D_1\subseteq D(h)$. In fact, we have that $b-a+1\leq \alpha-a=a$, so for each $j$ it holds $\{\xi,\xi+j\}\in D_j^{[n-h]}\subseteq D(h)$. And, since $n\geq 2\alpha+3$ and $h=2b-\alpha+2$, we have $n-h\geq 3\alpha-2b+1>2\alpha-2b$. So, $2\alpha-2b<\alpha<n-h$, thus for fixed $j$, the vertices $\{\xi,\xi+j\}$ are different. In consequence, $|D_1|=(b-a+1)(2\alpha-2b)=h(\alpha-b)$.
            
            On the other hand, $D_2\cap D(h)=\emptyset$ since no vertex in $D(h)$ has an element in $[n-h]$ and other in $[n-h+1..n]$. Thus, the vertices in $D_2$ are different and it turns out that $|D_2|=2(b-a+1)(2\alpha-2b)=2h(\alpha-b)$.

            If $\alpha$ is odd, let us consider the set $D=\left[D(h)\setminus D_1\right]\cup D_2$, where
            \begin{align*}
                D_1 =& \Bigg(\smashoperator[r]{\bigcup_{\substack{1\leq j\leq b-a\\1\leq\xi\leq 2\alpha - 2b}}}\ \big\{\{\xi,\xi+j\}\big\}\Bigg) \cup \Bigg( \smashoperator[r]{\bigcup_{1\leq\xi\leq \alpha - b} }\ \ \left\{\left\{\xi,\xi+\frac{n-h}{2}\right\}\right\}\Bigg),\\
                D_2 =& \Bigg(\smashoperator[r]{\bigcup_{\substack{1\leq j \leq b-a\\1\leq\xi\leq2\alpha-2b}}}\ \big\{\{\xi,n-h+2j-1\},\{\xi+j,n-h+2j\}\big\}\Bigg)\ \cup \\
                &\Bigg( \smashoperator[r]{\bigcup_{1\leq\xi\leq \alpha - b} }\ \ \left\{\left\{\xi,n\right\},\left\{\xi+\frac{n-h}{2},n\right\}\right\}\Bigg),
            \end{align*}
            and the sums in $D_1$ are taken modulo $n-h$.

            Let us see that $D_1\subseteq D(h)$. In fact, we have that $b-a\leq \alpha-1-a=a$, so for each $j$ it holds $\{\xi,\xi+j\}\in D_j^{[n-h]}\subseteq D(h)$. And, since $n\geq 2\alpha+3$ and $h=2b-\alpha+2$, we have $n-h\geq 3\alpha-2b+1>2\alpha-2b$. So, $2\alpha-2b<\alpha<n-h$, thus for fixed $j$, the vertices $\{\xi,\xi+j\}$ are different. On the other hand, for $1\leq\xi\leq \alpha-b<\frac{n-h}{2}$, we have $\left\{\xi,\xi+\frac{n-h}{2}\right\}\in \widehat{D}\subseteq D(h)$. In consequence, 
            $$|D_1|=(b-a)(2\alpha-2b)+(\alpha-b)=\underbrace{(2b-2a+1)}_{h}(\alpha-b)=h(\alpha-b).$$
            
            Besides, $D_2\cap D(h)=\emptyset$ since no vertex in $D(h)$ has an element in $[n-h]$ and other in $[n-h+1..n]$. Thus, the vertices in $D_2$ are different and it turns out that $|D_2|=2(b-a)(2\alpha-2b)+2(\alpha-b)=2h(\alpha-b)$.

            Either if $\alpha$ is even or odd, the set $D$ is a $k$-tuple dominating set with cardinality $k+2\alpha+1$. In fact, note that for each $x\in [n]$ we have
            $$i_x(D)=i_x(D(h))-i_x(D_1)+i_x(D_2).$$
            For $x\in [n-h]$, we have $i_x(D_1)=i_x(D_2)$, and in consequence 
            $$i_x(D)=i_x(D(h))=\alpha.$$
            And for $x\in [n-h+1..n]$, we have $i_x(D_1)=0$ and $i_x(D_2)=2\alpha-2b$, and in consequence 
            $$i_x(D)=\underbrace{i_x(D(h))}_{h-1}+2\alpha-2b=\alpha+1.$$
            We have $\alpha\leq i_x(D)\leq \alpha+1$ for each $x\in [n]$, and $\{x,y\}\in D$ for every pair $x,y\in X_{\alpha+1}$. Thus, $D$ is a $k$-tuple dominating set. Furthermore,
            $$|D|=|D(h)|-|D_1|+|D_2|=k+2\alpha+1.$$

            Therefore, either if $\alpha$ is even or odd, we obtain a $k$-tuple dominating set with cardinality $k+2\alpha+1$.
        \end{description}
    \end{proof}
\end{lemma}

\begin{example}
    \begin{itemize}
        \item Let $n=13$ and $k=15$. We have $n=\left\lceil\frac{2}{\alpha}k\right\rceil+3$ for $\alpha=3$, and $n\geq 9=2\alpha+3$. We are in the hypothesis of Proposition \ref{lem: kupla set cardinality k+2 alfa +1}. Thus, $\Kneser{13}{2}$ admits a $k$-tuple dominating set with cardinality $k+2\alpha+1=22$. Note that $k\equiv 0$ modulo $\alpha$, so it is enough to consider the set
        $$D=D^{8,3}\cup\binom{[9..13]}{2}.$$
        
        \item Let $n=18$ and $k=29$. We have $n=\left\lceil\frac{2}{\alpha}k\right\rceil+3$ for $\alpha=4$, and $n\geq 11=2\alpha+3$. We are in the hypothesis of Proposition \ref{lem: kupla set cardinality k+2 alfa +1}. Thus, $\Kneser{18}{2}$ admits a $k$-tuple dominating set with cardinality $k+2\alpha+1=38$. Note that $k\equiv 1$ modulo $\alpha$, so let us consider the set 
        $$D^*=D^{14,4}\cup \binom{[15..18]}{2},$$
        with $D^{14,4}$ as in example \ref{ej: ejemplos D(m,alfa)}. We obtain a $k$-tuple dominating set of $\Kneser{18}{2}$ deleting from $D^*$ the vertices
        $$\{1,2\}, \{2,3\}, \{1,3\}, \{2,4\}$$
        and adding the vertices
        $$\{1,15\},\{2,16\},\{2,15\},\{3,16\},\{1,17\},\{3,18\},\{2,17\},\{4,18\}.$$        
    \end{itemize}
\end{example}

Combining Proposition \ref{prop: lower bounds kdom(n,2) V2} with lemmas \ref{lem: kupla set cardinality k+2 alfa} and \ref{lem: kupla set cardinality k+2 alfa +1} yields the following result.

\begin{theorem}[\textbf{A}]\label{thm: kdom(n,2) con n en funcion de alfa}
    Let $\alpha,n\in \mathbb{N}$ with $\alpha\geq 2$ and $n\geq 2\alpha +3+(\alpha \mod 2)$. We have
    \begin{enumerate}
        \item $\kdom{n,2}{k} = k+2\alpha$, if $\frac{2}{\alpha}k+4\leq n < \frac{2}{\alpha-1}k+3$,
        \item $\kdom{n,2}{k} = k+2\alpha+1$, if $n=\left\lceil\frac{2}{\alpha}k\right\rceil+3$.
    \end{enumerate}
\end{theorem}


Let us see two remarks about the statement of Theorem \ref{thm: kdom(n,2) con n en funcion de alfa} (A).

Firstly, observe that $\alpha=\left\lceil\frac{2k}{n-3}\right\rceil$ if and only if $\frac{2}{\alpha}k+3\leq n < \frac{2}{\alpha-1}k+3$.
Thus, if $\alpha=\left\lceil\frac{2k}{n-3}\right\rceil$, then either $\frac{2}{\alpha}k+4 \leq n < \frac{2}{\alpha-1}k+3$, or $\frac{2}{\alpha}k+3\leq n < \frac{2}{\alpha}k+4$ (i.e. $n=\left\lceil\frac{2}{\alpha}k\right\rceil$+3).

Secondly, consider $k,n\in \mathbb{N}$ with $n\geq 7$ and $\frac{n-3}{2}< k\leq a\, \frac{n-3}{4}$ with $a=n-4$ if $n$ is even and $a=n-6+(n \mod 4)$ if $n$ is odd, and $\alpha=\left\lceil\frac{2k}{n-3}\right\rceil$. Let us see that $\alpha\geq 2$ and $n\geq 2\alpha+3+(\alpha\mod 2)$. The first inequality arises from the fact that $k>\frac{n-3}{2}$. In consequence $\alpha=\left\lceil\frac{2k}{n-3}\right\rceil\geq 2$. On the other hand, we have $k\leq a\,\frac{n-3}{4}$, which leads us to $\frac{2k}{n-3}\leq \frac{a}{2}$. Whether $n$ is even or odd, $a$ turns out to be even. Thus, $\frac{a}{2}\in \mathbb{N}$ and $\alpha=\left\lceil\frac{2k}{n-3}\right\rceil\leq \frac{a}{2}$, or equivalently $2\alpha\leq a$.
If $n$ is even, $a = n-4$, and it holds $n\geq 2\alpha+4$. 
If $n$ is odd, then $a=n-6+(n\mod 4)$, and it holds $n\geq 2\alpha+5$ if $n\equiv 1 \mod 4$, and $n\geq 2\alpha+3$ if $n\equiv 3 \mod 4$. Moreover, if $\alpha$ is odd and $n\equiv 3 \mod 4$, then $\frac{a}{2}=\frac{n-3}{2}$ is even. Thus, $\alpha\leq\frac{a}{2}-1$, and we have $n\geq 2\alpha+5$.
In any case, $n \geq 2\alpha+3+(\alpha \mod 2)$ as claimed.

Thus, Theorem \ref{thm: kdom(n,2) con n en funcion de alfa} (A) can be restated as follows.
\setcounter{theorem}{28}
\begin{theorem}[\textbf{B}]
\label{thm: alfas}
    Let $k,n\in \mathbb{N}$ with $n\geq 7$ and $\frac{n-3}{2}< k\leq a\, \frac{n-3}{4}$ where $a=n-4$ if $n$ is even and $a=n-6+(n \mod 4)$ if $n$ is odd. Let $\alpha=\left\lceil\frac{2k}{n-3}\right\rceil$. It holds
    \begin{enumerate}
        \item $\kdom{n,2}{k}=k+2\alpha$, if $\frac{2}{\alpha}k+4\leq n<\frac{2}{\alpha-1}k+3$;
        \item $\kdom{n,2}{k}=k+2\alpha+1$, if $n=\left\lceil\frac{2}{\alpha}k\right\rceil+3$.
    \end{enumerate}
\end{theorem}

It is worth noting that whereas in Theorem \ref{thm: alfas} (B) the values of $\kdom{n,2}{k}$ are computed for $\frac{n-3}{2}<k\leq a\,\frac{n-3}{4}$, with $a=n-4$ if $n$ is even and $a=n-6+(n\mod 4)$ if $n$ is odd, theorems \ref{thm: kdom = k+r} and \ref{thm: kdom(n,2)=k+3} provide the values of $\kdom{n,2}{k}$ for $k\leq \frac{n-3}{2}$. In consequence, we have determined $\kdom{n,2}{k}$ for all $2\leq k\leq a\, \frac{n-3}{4}$.

\section{Discussion}

Regarding $\gamma_{\times k}$-sets of $\Kneser{n}{2}$ for large values of $k$, recall that for $k=\binom{n-2}{2}+1$, we have $\kdom{n,2}{k}=\binom{n}{2}$, whereas for $k=\binom{n-2}{2}$, $\kdom{n,2}{k}=\binom{n}{2}-1$. Together with Theorem \ref{thm: large k} for $r=2$, we have the following. 

\begin{corollary}\label{coro: kdom(n,2) k grandes}
For $t\in \mathbb{N}$, if $k=\binom{n-2}{2}-t$ and $n\geq 2(t+3)-\left\lceil\frac{t+2}{2}\right\rceil$, then $\kdom{n,2}{k}=k+2n-4$.
\end{corollary}

Besides, in order to complete the study of $\gamma_{\times k}$-sets of $\Kneser{n}{2}$ for large values of $k$, we have also studied properties of the $\gamma_{\times k}$-sets of $\Kneser{n}{2}$ for $k\in\left\{\binom{n-4}{2}+2,\binom{n-3}{2}+1\right\}$. In fact, we obtain the following result. We omit the proof here for the sake of readability (it can be found in \cite{ktuple2023v1}).

\begin{proposition}\label{prop:partcases}
    \begin{itemize}
        \item  If $k=\binom{n-4}{2}+2$ and $6\leq n\leq 10$, then $\kdom{n,2}{k}=\binom{n-2}{2}+1$. Moreover, $D=\binom{[n-2]}{2}\cup\big\{\{n-1,n\}\big\}$ is a $\gamma_{\times k}$-set.
        \item  If $k=\binom{n-3}{2}+1$ and $n\geq 5$, then $\kdom{n,2}{k}=\binom{n-1}{2}$. Moreover, $D=\binom{[n-1]}{2}$ is a $\gamma_{\times k}$-set.
    \end{itemize}
\end{proposition}

Notice that although the set $\binom{[n-2]}{2}\cup\big\{\{n-1,n\}\big\}$ given in the first item of Proposition \ref{prop:partcases} is a  $k$-tuple dominating set of $\Kneser{n}{2}$ for every $n\geq 6$, it is not true that it is a $\gamma_{\times k}$-set for every $n$. For instance, for $n=11$ and $k=23$ the set
$$D=\binom{[6]}{2}\cup \binom{[7..11]}{2}\cup \tilde{D},$$
where
$$\tilde{D}=\big\{\{1,7\},\{2,7\},\{2,8\},\{3,8\},\{3,9\},\{4,9\},\{4,10\},\{5,10\},\{5,11\}\big\},$$
is a $k$-tuple dominating set with cardinality $\binom{n-2}{2}=\binom{9}{2}=36$.

With the aim of summarizing the results in this paper for $\Kneser{n}{2}$, in Appendix \ref{sec:table} we give a table that contains several values for $\kdom{n,2}{k}$.
Some of these values are obtained taking into account Lemma \ref{lemma: cotas indices}, Theorem \ref{thm: monotonia respecto a n}, and upper bounds obtained by Integer Linear Programming. In this regard, we consider the following ILP formulation for the $k$-tuple dominating problem. The set of variables is given by $\{x_u, u\in V\}$, where
$$x_u=\begin{cases}
    1, &\text{if }u\in D,\\
    0, &\text{otherwise,}
\end{cases}$$
and the ILP model is formulated as
\begin{equation}\label{eq: ILP}
    \begin{aligned}
        \text{min} & \sum_{u\in V} x_u\\
        \text{s/t} & \sum_{u\in \neighc{v}{}}x_u\geq k, && \forall\, v\in V,\\
        & x_v\in\{0,1\}, && \forall\, v\in V.
    \end{aligned}
\end{equation}
We implemented it in CPLEX solver \cite{CPLEX} to obtain $k$-tuple dominating sets that allowed us to determine an upper bound of $\kdom{n,2}{k}$ for certain values of $n$ and $k$, which together with the lower bounds obtained using  Lemma \ref{lemma: cotas indices} and Theorem \ref{thm: monotonia respecto a n}, were tight. Notice that we have determined $\kdom{n,2}{k}$ for every $n$ and $k\leq18$. For the general case, fixed $k$, the amount of values of $n$ for which $\kdom{n,2}{k}$ remains unknown is $\Theta(\sqrt{k})$. 
However, for some of them we have upper bounds arising from solving the ILP.







\appendix

\section{\texorpdfstring{$\gamma_{\times k}$-sets of the odd graph $\Kneser{11}{5}$}{}}\label{sec: K(11,5)}

Consider the Kneser graph $\Kneser{11}{5}$. We will see that for any $k\in [7]$, we have 
\begin{equation}
    \label{eq: kdom(11,5)}
    \kdom{11,5}{k}=k\packing{11,5}=66k.
\end{equation}
For $k = 1$, an minimum dominating set is also a maximum $2$-packing, which is given by a Steiner System $S(4,5,11)$.

We obtain a Steiner system $S=S(4,5,11)$ with base set $[n]$ following the construction given in \cite{grannell1994introduction}. Let us consider two sets $A=\{1,2,3,4,5\}$, $B=\{6,7,8,9,10,11\}$. Suppose that $A$ is itself a block, and consider the following scheme.
\begin{center}
\begin{tabular}{c c c c}
    1 & (6,7) & (8,9) & (10,11) \\
    2 & (6,8) & (9,10) & (7,11) \\
    3 & (6,9) & (7,10) & (8,11) \\
    4 & (6,10) & (7,8) & (9,11) \\
    5 & (6,11) & (8,10) & (7,9)
\end{tabular}
\end{center}
In the set $S$ there are $30$ blocks of the form $aaabb$, $20$ blocks of the form $aabbb$, $15$ of the form $abbbb$, and one block $aaaaa$, where $a$ denotes an element in $A$ and $b$ denotes an element in $B$. 
Notice that in each row of the scheme there is an element $a\in A$ and three disjoint pairs of elements of $B$. 
We form three blocks from each row taking the element $a$ and two of the three pairs. Doing so, we get the $15$ blocks of the form $abbbb$.
Now, for any three elements $bbb$ of $B$, taking into account that any four elements are contained in exactly one block of $S$, there are exactly two elements of $A$ for which $aabbb$ is a block. For instance, for $\{6,7,8\}$, we have that $\{1,6,7,8,9\}$, $\{2,6,7,8,11\}$ and $\{4,6,7,8,10\}$ are blocks in $S$ of the form $abbbb$. Thus, the block in $S$ of the form $aabbb$ containing $\{6,7,8\}$ is $\{3,5,6,7,8\}$.
Finally, for any two elements $bb$ of $B$, with a similar reasoning we get two blocks of the form $aaabb$. For example, for the set $\{6,7\}$ we have $\{3,5,6,7,8\}$, $\{2,4,6,7,9\}$, $\{2,5,6,7,10\}$ and $\{3,4,6,7,11\}$ are blocks of the form $aabbb$ which are in $S$. Thus, the blocks in $S$ of the form $aaabb$ containing the pair $\{6,7\}$ are $\{1,2,3,6,7\}$ and $\{1,4,5,6,7\}$.

With this procedure, we get a Steiner system $S(4,5,11)$. $S$ turns out to be a $\gamma_{\times 1}$-set of the Kneser graph $\Kneser{11}{5}$ as well as an maximum $2$-packing.

For $k=2$, we obtain a $\gamma_{\times 2}$-set of $\Kneser{11}{5}$ by taking $S_1\cup S_2$ where $S_1$ and $S_2$ and two disjoint Steiner systems $S(4,5,11)$. It is enough to consider the construction above with the following schemes.
\begin{multicols}{2}
\begin{center}
Scheme for $S_1$

\begin{tabular}{c c c c}
    1 & (6,7) & (8,9) & (10,11) \\
    2 & (6,8) & (9,10) & (7,11) \\
    3 & (6,9) & (7,10) & (8,11) \\
    4 & (6,10) & (7,8) & (9,11) \\
    5 & (6,11) & (8,10) & (7,9)
\end{tabular}

Scheme for $S_2$

\begin{tabular}{c c c c}
    6 & (1,5) & (2,4) & (3,7) \\
    8 & (1,2) & (3,5) & (4,7) \\
    9 & (1,4) & (2,3) & (5,7) \\
    10 & (1,3) & (2,7) & (4,5) \\
    11 & (1,7) & (2,5) & (3,4) 
\end{tabular}
\end{center}
\end{multicols}

For $k=3$, we found a $k$-tuple dominating set $D$ with cardinality $198=3\cdot 66$. Below, we list the vertices in $D$.
\begin{multicols}{5}
\begin{center}
    \{1, 2, 3, 4, 5\} \\
    \{1, 2, 3, 4, 7\} \\
    \{1, 2, 3, 4, 8\} \\
    \{1, 2, 3, 5, 9\} \\
    \{1, 2, 3, 5, 11\} \\
    \{1, 2, 3, 6, 7\} \\
    \{1, 2, 3, 6, 8\} \\
    \{1, 2, 3, 6, 11\} \\
    \{1, 2, 3, 7, 10\} \\
    \{1, 2, 3, 8, 9\} \\
    \{1, 2, 3, 9, 10\} \\
    \{1, 2, 3, 10, 11\} \\
    \{1, 2, 4, 5, 6\} \\
    \{1, 2, 4, 5, 9\} \\
    \{1, 2, 4, 6, 8\} \\
    \{1, 2, 4, 6, 10\} \\
    \{1, 2, 4, 7, 9\} \\
    \{1, 2, 4, 7, 11\} \\
    \{1, 2, 4, 8, 10\} \\
    \{1, 2, 4, 9, 11\} \\
    \{1, 2, 4, 10, 11\} \\
    \{1, 2, 5, 6, 8\} \\
    \{1, 2, 5, 6, 10\} \\
    \{1, 2, 5, 7, 8\} \\
    \{1, 2, 5, 7, 9\} \\
    \{1, 2, 5, 7, 10\} \\
    \{1, 2, 5, 8, 11\} \\
    \{1, 2, 5, 10, 11\} \\
    \{1, 2, 6, 7, 9\} \\
    \{1, 2, 6, 7, 11\} \\
    \{1, 2, 6, 9, 10\} \\
    \{1, 2, 6, 9, 11\} \\
    \{1, 2, 7, 8, 10\} \\
    \{1, 2, 7, 8, 11\} \\
    \{1, 2, 8, 9, 10\} \\
    \{1, 2, 8, 9, 11\} \\
    \{1, 3, 4, 5, 6\} \\
    \{1, 3, 4, 5, 10\} \\
    \{1, 3, 4, 6, 7\} \\
    \{1, 3, 4, 6, 9\} \\
    \{1, 3, 4, 7, 11\} \\
    \{1, 3, 4, 8, 9\} \\
    \{1, 3, 4, 8, 10\} \\
    \{1, 3, 4, 9, 11\} \\
    \{1, 3, 4, 10, 11\} \\
    \{1, 3, 5, 6, 7\} \\
    \{1, 3, 5, 6, 10\} \\
    \{1, 3, 5, 7, 8\} \\
    \{1, 3, 5, 7, 9\} \\
    \{1, 3, 5, 8, 10\} \\
    \{1, 3, 5, 8, 11\} \\
    \{1, 3, 5, 9, 11\} \\
    \{1, 3, 6, 8, 9\} \\
    \{1, 3, 6, 8, 11\} \\
    \{1, 3, 6, 9, 10\} \\
    \{1, 3, 6, 10, 11\} \\
    \{1, 3, 7, 8, 10\} \\
    \{1, 3, 7, 8, 11\} \\
    \{1, 3, 7, 9, 10\} \\
    \{1, 3, 7, 9, 11\} \\
    \{1, 4, 5, 6, 11\} \\
    \{1, 4, 5, 7, 8\} \\
    \{1, 4, 5, 7, 10\} \\
    \{1, 4, 5, 7, 11\} \\
    \{1, 4, 5, 8, 9\} \\
    \{1, 4, 5, 8, 11\} \\
    \{1, 4, 5, 9, 10\} \\
    \{1, 4, 6, 7, 8\} \\
    \{1, 4, 6, 7, 10\} \\
    \{1, 4, 6, 8, 11\} \\
    \{1, 4, 6, 9, 10\} \\
    \{1, 4, 6, 9, 11\} \\
    \{1, 4, 7, 8, 9\} \\
    \{1, 4, 7, 9, 10\} \\
    \{1, 4, 8, 10, 11\} \\
    \{1, 5, 6, 7, 9\} \\
    \{1, 5, 6, 7, 11\} \\
    \{1, 5, 6, 8, 9\} \\
    \{1, 5, 6, 8, 10\} \\
    \{1, 5, 6, 9, 11\} \\
    \{1, 5, 7, 10, 11\} \\
    \{1, 5, 8, 9, 10\} \\
    \{1, 5, 9, 10, 11\} \\
    \{1, 6, 7, 8, 9\} \\
    \{1, 6, 7, 8, 10\} \\
    \{1, 6, 7, 10, 11\} \\
    \{1, 6, 8, 10, 11\} \\
    \{1, 7, 8, 9, 11\} \\
    \{1, 7, 9, 10, 11\} \\
    \{1, 8, 9, 10, 11\} \\
    \{2, 3, 4, 5, 7\} \\
    \{2, 3, 4, 5, 8\} \\
    \{2, 3, 4, 6, 9\} \\
    \{2, 3, 4, 6, 10\} \\
    \{2, 3, 4, 6, 11\} \\
    \{2, 3, 4, 7, 10\} \\
    \{2, 3, 4, 8, 11\} \\
    \{2, 3, 4, 9, 10\} \\
    \{2, 3, 4, 9, 11\} \\
    \{2, 3, 5, 6, 9\} \\
    \{2, 3, 5, 6, 10\} \\
    \{2, 3, 5, 6, 11\} \\
    \{2, 3, 5, 7, 10\} \\
    \{2, 3, 5, 7, 11\} \\
    \{2, 3, 5, 8, 9\} \\
    \{2, 3, 5, 8, 10\} \\
    \{2, 3, 6, 7, 8\} \\
    \{2, 3, 6, 7, 9\} \\
    \{2, 3, 6, 8, 10\} \\
    \{2, 3, 7, 8, 9\} \\
    \{2, 3, 7, 8, 11\} \\
    \{2, 3, 7, 9, 11\} \\
    \{2, 3, 8, 10, 11\} \\
    \{2, 3, 9, 10, 11\} \\
    \{2, 4, 5, 6, 7\} \\
    \{2, 4, 5, 6, 11\} \\
    \{2, 4, 5, 7, 8\} \\
    \{2, 4, 5, 8, 10\} \\
    \{2, 4, 5, 9, 10\} \\
    \{2, 4, 5, 9, 11\} \\
    \{2, 4, 5, 10, 11\} \\
    \{2, 4, 6, 7, 9\} \\
    \{2, 4, 6, 7, 10\} \\
    \{2, 4, 6, 8, 9\} \\
    \{2, 4, 6, 8, 11\} \\
    \{2, 4, 7, 8, 9\} \\
    \{2, 4, 7, 8, 11\} \\
    \{2, 4, 7, 10, 11\} \\
    \{2, 4, 8, 9, 10\} \\
    \{2, 5, 6, 7, 8\} \\
    \{2, 5, 6, 7, 11\} \\
    \{2, 5, 6, 8, 9\} \\
    \{2, 5, 6, 9, 10\} \\
    \{2, 5, 7, 9, 10\} \\
    \{2, 5, 7, 9, 11\} \\
    \{2, 5, 8, 9, 11\} \\
    \{2, 5, 8, 10, 11\} \\
    \{2, 6, 7, 8, 10\} \\
    \{2, 6, 7, 10, 11\} \\
    \{2, 6, 8, 9, 11\} \\
    \{2, 6, 8, 10, 11\} \\
    \{2, 6, 9, 10, 11\} \\
    \{2, 7, 8, 9, 10\} \\
    \{2, 7, 9, 10, 11\} \\
    \{3, 4, 5, 6, 8\} \\
    \{3, 4, 5, 6, 9\} \\
    \{3, 4, 5, 7, 9\} \\
    \{3, 4, 5, 7, 11\} \\
    \{3, 4, 5, 8, 11\} \\
    \{3, 4, 5, 9, 10\} \\
    \{3, 4, 5, 10, 11\} \\
    \{3, 4, 6, 7, 8\} \\
    \{3, 4, 6, 7, 11\} \\
    \{3, 4, 6, 8, 10\} \\
    \{3, 4, 6, 10, 11\} \\
    \{3, 4, 7, 8, 9\} \\
    \{3, 4, 7, 8, 10\} \\
    \{3, 4, 7, 9, 10\} \\
    \{3, 4, 8, 9, 11\} \\
    \{3, 5, 6, 7, 8\} \\
    \{3, 5, 6, 7, 10\} \\
    \{3, 5, 6, 8, 11\} \\
    \{3, 5, 6, 9, 11\} \\
    \{3, 5, 7, 8, 9\} \\
    \{3, 5, 7, 10, 11\} \\
    \{3, 5, 8, 9, 10\} \\
    \{3, 5, 9, 10, 11\} \\
    \{3, 6, 7, 9, 10\} \\
    \{3, 6, 7, 9, 11\} \\
    \{3, 6, 7, 10, 11\} \\
    \{3, 6, 8, 9, 10\} \\
    \{3, 6, 8, 9, 11\} \\
    \{3, 7, 8, 10, 11\} \\
    \{3, 8, 9, 10, 11\} \\
    \{4, 5, 6, 7, 9\} \\
    \{4, 5, 6, 7, 10\} \\
    \{4, 5, 6, 8, 9\} \\
    \{4, 5, 6, 8, 10\} \\
    \{4, 5, 6, 10, 11\} \\
    \{4, 5, 7, 8, 10\} \\
    \{4, 5, 7, 9, 11\} \\
    \{4, 5, 8, 9, 11\} \\
    \{4, 6, 7, 8, 11\} \\
    \{4, 6, 7, 9, 11\} \\
    \{4, 6, 8, 9, 10\} \\
    \{4, 6, 9, 10, 11\} \\
    \{4, 7, 8, 10, 11\} \\
    \{4, 7, 9, 10, 11\} \\
    \{4, 8, 9, 10, 11\} \\
    \{5, 6, 7, 8, 11\} \\
    \{5, 6, 7, 9, 10\} \\
    \{5, 6, 8, 10, 11\} \\
    \{5, 6, 9, 10, 11\} \\
    \{5, 7, 8, 9, 10\} \\
    \{5, 7, 8, 9, 11\} \\
    \{5, 7, 8, 10, 11\} \\
    \{6, 7, 8, 9, 10\} \\
    \{6, 7, 8, 9, 11\}
\end{center}
\end{multicols}

For $k\in\{1,2,3\}$ we gave a $\gamma_{\times k}$-set $D_k$ with cardinality $66k$. By Remark \ref{rem: complemento de kupla dom cota packing}, we have that $D_{7-k}=V\setminus D_k$ is a $(7-k)$-tuple dominating set with cardinality $|D_{7-k}|=|V|-66k=66(7-k)$. Finally, the set $D_7=V$ is a $\gamma_{\times 7}$-set with cardinality $66\cdot 7$. Thus, (\ref{eq: kdom(11,5)}) holds for every $k\in [7]$.

\section{\texorpdfstring{Table for $\kdom{n,2}{k}$}{Table for kdom(n,2)}}\label{sec:table}

Table \ref{tab: tabla K(n,2)} shows the values of $\kdom{n,2}{k}$ stated by the results in this paper. Rows correspond to values of $k$ while columns correspond to values of $n$. The superscript in each entry of the table indicates from which result it follows, and the bounds arise from $k$-tuple dominating sets obtained using CPLEX solver for (\ref{eq: ILP}) and monotonicity.
\begin{enumerate}[($a$)]
    \item Domination number \cite{ivanvco1993domination}.
    \item Theorem \ref{thm: alfas} (B).1.
    \item Theorem \ref{thm: alfas} (B).2.
    \item Corollary \ref{coro: kdom(n,2) k grandes}.
    \item Proposition \ref{prop:partcases}.
    \item Lower bound from Lemma \ref{lemma: cotas indices} and upper bound from $\gamma_{\times k}$-sets found by ILP (\ref{eq: ILP}).
    \item Lower bound from Monotonicity (Theorem \ref{thm: monotonia respecto a n}) and upper bound from $\gamma_{\times k}$-sets found by ILP (\ref{eq: ILP}).    
\end{enumerate}

\begin{table}[H]
    \centering
    \includegraphics[height=0.92\textheight]{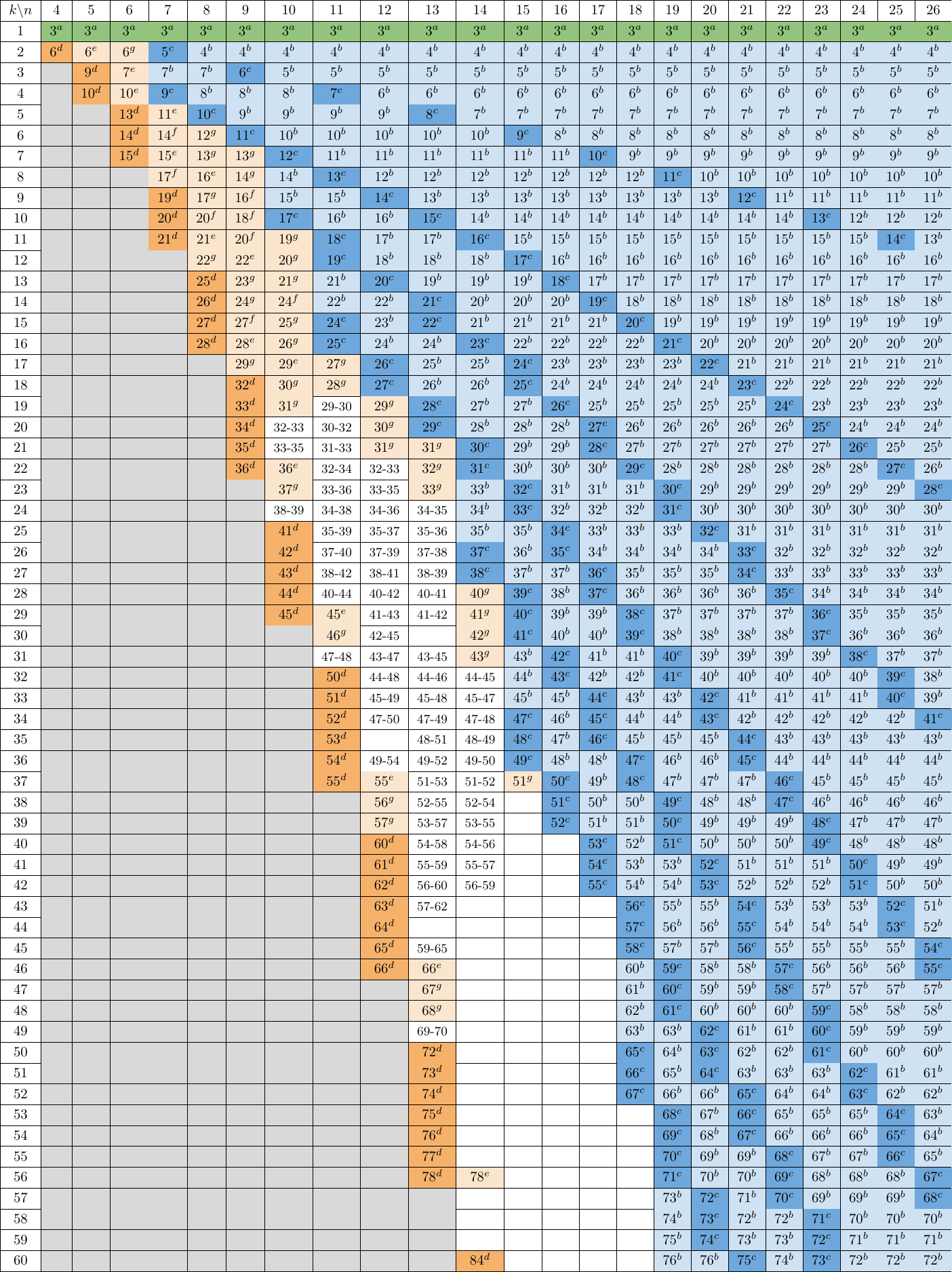}
    \caption{$\kdom{n,2}{k}$ for some values of $n$ and $k$.}
    \label{tab: tabla K(n,2)}
\end{table}

\newpage

\bibliographystyle{acm}
\bibliography{biblio}

\begin{thebibliography}{10}

\bibitem{barrau1908combinatory}
{\sc Barrau, J.}
\newblock On the combinatory problem of {Steiner}.
\newblock {\em Koninklijke Nederlandse Akademie van Wetenschappen Proceedings Series B Physical Sciences 11\/} (1908), 352--360.

\bibitem{bedo2023geodetic}
{\sc Bedo, M., Leite, J.~V., Oliveira, R.~A., and Protti, F.}
\newblock Geodetic convexity and {Kneser} graphs.
\newblock {\em Applied Mathematics and Computation 449\/} (2023), 127964.

\bibitem{biggs1979}
{\sc Biggs, N.}
\newblock Some odd graph theory.
\newblock {\em Annals of the New York Academy of Sciences 319}, 1 (1979), 71--81.

\bibitem{2packing-BKPG-2022}
{\sc Bozovic, D., Kelenc, A., Peterin, I., and González~Yero, I.}
\newblock Incidence dimension and $2$-packing number in graphs.
\newblock {\em RAIRO - Operations Research 56}, 1 (2022), 199--211.

\bibitem{bresar2019grundy}
{\sc Bre\v{s}ar, B., Kos, T., and Torres, P.~D.}
\newblock Grundy domination and zero forcing in {Kneser} graphs.
\newblock {\em ARS MATHEMATICA CONTEMPORANEA 17\/} (2019), 419--430.

\bibitem{cornet2023kdom}
{\sc Bre{š}ar, B., Cornet, M.~G., Dravec, T., and Henning, M.~A.}
\newblock $k$-domination invariants on {Kneser} graphs.
\newblock {\em arXiv e-prints\/} (2023), arxiv:2312.15464.

\bibitem{martinez2022note}
{\sc Cabrera~Mart{\'\i}nez, A.}
\newblock A note on the $k$-tuple domination number of graphs.
\newblock {\em ARS MATHEMATICA CONTEMPORANEA 5\/} (2022).

\bibitem{2packingCKMR-2011}
{\sc Castro, A., Klavžar, S., Mollard, M., and Rho, Y.}
\newblock On the domination number and the $2$-packing number of {F}ibonacci cubes and {L}ucas cubes.
\newblock {\em Computers \& Mathematics with Applications 61}, 9 (2011), 2655--2660.

\bibitem{chang2008upper}
{\sc Chang, G.~J.}
\newblock The upper bound on $k$-tuple domination numbers of graphs.
\newblock {\em European Journal of Combinatorics 29}, 5 (2008), 1333--1336.

\bibitem{2packingCG-2017}
{\sc Clarke, N.~E., and Gallant, R.~P.}
\newblock On $2$-limited packings of complete grid graphs.
\newblock {\em Discrete Mathematics 340}, 7 (2017), 1705–1715.

\bibitem{colbourn2006steiner}
{\sc Colbourn, C.~J., and Mathon, R.}
\newblock Steiner systems.
\newblock In {\em Handbook of Combinatorial Designs}. Chapman and Hall/CRC, 2006, pp.~128--135.

\bibitem{ktuple2023v1}
{\sc Cornet, M.~G., and Torres, P.}
\newblock $k$-tuple domination on {Kneser} graphs.
\newblock {\em arXiv e-prints\/} (2023), arXiv:2308.15603v1.

\bibitem{ekr-61}
{\sc Erd\H{o}s, P., Ko, C., and Rado, R.}
\newblock Intersection theorem for system of finite sets.
\newblock {\em Quart.\ J.\ Math. 12\/} (1961), 313--318.

\bibitem{qAnalog}
{\sc Fallat, S., Joshi, N., Maleki, R., Meagher, K., Mojallal, S.~A., Nasserasr, S., Shirazi, M.~N., Razafimahatratra, A.~S., and Stevens, B.}
\newblock The $q$-analogue of zero forcing for certain families of graphs.
\newblock {\em Discrete Applied Mathematics 348\/} (2024), 91--104.

\bibitem{gagarin2008generalised}
{\sc Gagarin, A., and Zverovich, V.~E.}
\newblock A generalised upper bound for the $k$-tuple domination number.
\newblock {\em Discrete Mathematics 308}, 5-6 (2008), 880--885.

\bibitem{GenPosKneser}
{\sc Ghorbani, M., Klavžar, S., Maimani, H.~R., Momeni, M., Mahid, F.~R., and Rus, G.}
\newblock The general position problem on {Kneser} graphs and on some graph operations.
\newblock {\em Discussiones Mathematicae Graph Theory 41\/} (2021), 1199--1213.

\bibitem{gorodezky2007dominating}
{\sc Gorodezky, I.}
\newblock Dominating sets in {Kneser} graphs.
\newblock Master's thesis, University of Waterloo, 2007.

\bibitem{grannell1994introduction}
{\sc Grannell, M., and Griggs, T.}
\newblock An introduction to {S}teiner systems.
\newblock {\em Mathematical Spectrum 26}, 3 (1994), 74--80.

\bibitem{hammond1975perfect}
{\sc Hammond, P., and Smith, D.}
\newblock Perfect codes in the graphs {$O_k$}.
\newblock {\em Journal of Combinatorial Theory, Series B 19}, 3 (1975), 239--255.

\bibitem{hw-2003}
{\sc Hartman, C., and West, D.~B.}
\newblock Covering designs and domination in {Kneser} graphs.
\newblock {\em unpublished manuscript\/} (2003).

\bibitem{TreewidthKneser}
{\sc Harvey, D.~J., and Wood, D.~R.}
\newblock Treewidth of the {Kneser} graph and the {E}rd\"os-{K}o-{R}ado theorem.
\newblock {\em The electronic Journal of Combinatorics 21:1\/} (2014), \#P1.48.

\bibitem{haynes2020topics}
{\sc Haynes, T.~W., Hedetniemi, S.~T., and Henning, M.~A.}
\newblock {\em Topics in domination in graphs}, vol.~64.
\newblock Springer, 2020.

\bibitem{hhs-1998b}
{\sc Haynes, T.~W., Hedetniemi, S.~T., and Slater, P.~J.}, Eds.
\newblock {\em Domination in graphs}, vol.~209 of {\em Monographs and Textbooks in Pure and Applied Mathematics}.
\newblock Marcel Dekker, Inc., New York, 1998.
\newblock Advanced topics.

\bibitem{hhs-1998a}
{\sc Haynes, T.~W., Hedetniemi, S.~T., and Slater, P.~J.}
\newblock {\em Fundamentals of domination in graphs}, vol.~208 of {\em Monographs and Textbooks in Pure and Applied Mathematics}.
\newblock Marcel Dekker, Inc., New York, 1998.

\bibitem{CPLEX}
{\sc {IBM}}.
\newblock {ILOG CPLEX Optimization Studio 22.1.0}.
\newblock \url{https://www.ibm.com/products/ilog-cplex-optimization-studio}, 2022.
\newblock Accessed 2024-03-30.

\bibitem{ivanvco1993domination}
{\sc Ivan{\v{c}}o, J., and Zelinka, B.}
\newblock Domination in {Kneser} graphs.
\newblock {\em Mathematica Bohemica 118}, 2 (1993), 147--152.

\bibitem{jalilolghadr2023total}
{\sc Jalilolghadr, P., and Behtoei, A.}
\newblock Total dominator chromatic number of {K}neser graphs.
\newblock {\em AKCE International Journal of Graphs and Combinatorics 20}, 1 (2023), 52--56.

\bibitem{2packing-JSR-2012}
{\sc Junosza-Szaniawski, K., and Rzazewski, P.}
\newblock On the number of $2$-packings in a connected graph.
\newblock {\em Discrete Mathematics 312}, 23 (2012), 3444--3450.

\bibitem{kneser1955aufgabe}
{\sc Kneser, M.}
\newblock Aufgabe 300.
\newblock {\em Jahresbericht Deutschen Math. Verein 58}, 2 (1955).

\bibitem{liao2003k}
{\sc Liao, C.-S., and Chang, G.~J.}
\newblock $k$-tuple domination in graphs.
\newblock {\em Information Processing Letters 87\/} (2003), 45--50.

\bibitem{lov-78}
{\sc Lov\'{a}sz, L.}
\newblock Kneser's conjecture, chromatic numbers and homotopy.
\newblock {\em J. Combin. Theory (A) 25\/} (1978), 319--324.

\bibitem{mat-04}
{\sc Matou\v{s}ek, J.}
\newblock A combinatorial proof of {Kneser}'s conjecture.
\newblock {\em Combinatorica 24\/} (2004), 163--170.

\bibitem{MeirMoon1975}
{\sc Meir, A., and Moon, J.~W.}
\newblock Relations between packing and covering numbers of a tree.
\newblock {\em Pacific J. Math. 61}, 1 (1975), 225--233.

\bibitem{ostergaard2014bounds}
{\sc {\"O}sterg{\aa}rd, P.~R., Shao, Z., and Xu, X.}
\newblock Bounds on the domination number of {Kneser} graphs.
\newblock {\em ARS MATHEMATICA CONTEMPORANEA 9}, 2 (2014), 187--195.

\bibitem{przybylo2013upper}
{\sc Przyby{\l}o, J.}
\newblock On upper bounds for multiple domination numbers of graphs.
\newblock {\em Discrete Applied Mathematics 161}, 16-17 (2013), 2758--2763.

\bibitem{ramanan1997proof}
{\sc Ramanan, G.~V.}
\newblock Proof of a conjecture of {Frankl} and {F{\"u}redi}.
\newblock {\em journal of combinatorial theory, Series A 79}, 1 (1997), 53--67.

\bibitem{RAUTENBACH200798}
{\sc Rautenbach, D., and Volkmann, L.}
\newblock New bounds on the $k$-domination number and the $k$-tuple domination number.
\newblock {\em Applied Mathematics Letters 20}, 1 (2007), 98--102.

\bibitem{valencia2005diameter}
{\sc Valencia-Pabon, M., and Vera, J.-C.}
\newblock On the diameter of {Kneser} graphs.
\newblock {\em Discrete mathematics 305\/} (2005), 383--385.

\bibitem{west2001introduction}
{\sc West, D.~B., et~al.}
\newblock {\em Introduction to graph theory}, vol.~2.
\newblock Prentice hall Upper Saddle River, 2001.

\bibitem{zec2023several}
{\sc Zec, T., and Grbi{\'c}, M.}
\newblock Several {R}oman domination graph invariants on {K}neser graphs.
\newblock {\em Discrete Mathematics \& Theoretical Computer Science 25:1}, Graph Theory (2023).

\bibitem{zverovich2008k}
{\sc Zverovich, V.}
\newblock The $k$-tuple domination number revisited.
\newblock {\em Applied Mathematics Letters 21}, 10 (2008), 1005--1011.

\end{thebibliography}

\end{document}